\documentclass[12pt,a4paper,reqno]{amsart}

\usepackage[utf8]{inputenc}
\usepackage{amsmath, amsthm}
\usepackage{graphicx}
\usepackage{amssymb}
\usepackage{lmodern} 
\usepackage{amsfonts}
\usepackage{bbm,bm}
\usepackage{esint}
\PassOptionsToPackage{dvipsnames}{xcolor}
\usepackage{tikz}
\usepackage{amsthm}
\usepackage{epsfig}
\usepackage{floatflt}
\usepackage[ngerman,english]{babel}
\usepackage{comment}
\usepackage[T1]{fontenc}        
\usepackage{url}
\usepackage{amssymb}
\usepackage{color}
\usepackage{setspace}
\usepackage{nicefrac}
\usepackage[colorlinks, linkcolor = black, citecolor = black, filecolor = black, urlcolor = blue]{hyperref}
\usepackage{esint} 
\usepackage{fancyref}
\usepackage{ifthen}
\usepackage{xcolor}

\usepackage{commath}

\allowdisplaybreaks[2]

\begin{document}
	\numberwithin{equation}{section}
	\newcommand{\proofend}{\begin{flushright}$\Box$\end{flushright}}
	\renewcommand{\proof}{\textbf{Proof}: }
	\newcommand{\R}{\mathbb{R}}
	\newcommand{\C}{\mathbb{C}}
	\newcommand{\N}{\mathbb{N}}
	\newcommand{\Z}{\mathbb{Z}}
	\newcommand{\Suml}{\sum\limits}
	\newcommand{\Cupl}{\bigcup\limits}
	\newcommand{\Capl}{\bigcap\limits}
	\newcommand{\Intl}{\int\limits}
	\newcommand{\Liml}{\lim\limits}
	\newcommand{\supl}{\sup\limits}
	\newcommand{\infl}{\inf\limits}
	\newcommand{\rL}{\mathrm{L}}
	\newcommand{\cL}{\mathcal{L}}
	\newcommand{\cW}{\mathcal{W}}
	\newcommand{\cD}{\mathcal{D}}
	\newcommand{\sS}{\mathcal{S}}
	\newcommand{\cO}{\mathcal{O}}
	\newcommand{\cX}{\mathcal{X}}
	\newcommand{\cH}{\mathcal{H}}
	\newcommand{\cG}{\mathcal{G}}
	\newcommand{\BV}{\mathrm{BV}}
	\newcommand{\rP}{\mathrm{P}}
	\newcommand{\divv}{\mathrm{div}\,}
	\newcommand{\cE}{\mathcal{E}}
	\newcommand{\M}{\mathcal{M}}
	\newcommand{\B}{\mathcal{B}}
	\newcommand{\X}{\mathcal{X}}
	\newcommand{\F}{\mathcal{F}}
	\newcommand{\e}{\varepsilon}
	\newcommand{\s}{\sigma}
	\newcommand{\g}{\gamma}
	\newcommand{\ca}{\operatorname{cap}}
	\newcommand{\diam}{\operatorname{diam}}
	\newcommand{\maxi}{\operatorname{max}}
	\newcommand{\mini}{\operatorname{min}}
	\newcommand{\rad}{\operatorname{rad}}
	\newcommand{\dist}{\operatorname{dist}}
	
	\newcommand{\dx}{\dif x}
	\newcommand{\dy}{\dif y}
	\newcommand{\dt}{\dif t}
	\newcommand{\ds}{\dif s}
	\newcommand{\dbl}{\dif \bl}
	\newcommand{\dble}{\dif \ble}
	\newcommand{\dblt}{\dif \blt}
	\newcommand{\dblet}{\dif \blet}
	\newcommand{\intsp}[1]{\int_{\cO} #1 \dif x}
	\newcommand{\inttimet}[1]{\int_{0}^{t} #1 \, \dif s}
	\newcommand{\inttimeT}[1]{\int_{0}^{T} #1 \, \dif s}
	\newcommand{\inttimetra}[1]{\int_{0}^{t} #1 \, \dif \beta_{k}}
	\newcommand{\intsptit}[1]{\int_{0}^{t}\int_{\cO} #1 \, \dif x \dif s}
	\newcommand{\intsptist}[1]{\int_{s}^{t}\int_{\cO} #1 \, \dif x \dif s}
	\newcommand{\intsptiT}[1]{\int_{0}^{T}\int_{\cO} #1 \, \dif x \dif s}
	\newcommand{\intsptitra}[1]{\int_{0}^{t}\int_{\cO} #1 \, \dif x \dif \beta_{k}^{\e}}
	\newcommand{\intsptiTra}[1]{\int_{0}^{T}\int_{\cO} #1 \, \dif x \dif \beta_{k}^{\e}}
	
	\newcommand{\intEwert}[1]{\mathbb{E} \left[ \int_{0}^{t}\int_{\cO} #1 \, \dif x \dif s \right]}
	
	\newcommand{\betr}[1]{\left| #1 \right|}
	\newcommand{\lpr}[2]{L^{#1}(#2)}
	\newcommand{\Hsob}[2]{H^{#1}(#2)}
	\newcommand{\Hsobper}[2]{H_{\text{per}}^{#1}(#2)}
	\newcommand{\sob}[2]{W^{#1}(#2)}
	\newcommand{\li}[2]{\underset{#1 \rightarrow #2 }{\lim}}
	\newcommand{\Ew}[1]{\mathbb{E}\left[  #1  \right]}
	\newcommand{\Ewt}[1]{\ti{\mathbb{E}}\left[  #1  \right]}
	\newcommand{\prob}[1]{\mathbb{P}\left[  #1  \right]}
	\newcommand{\ti}[1]{\tilde{#1}}
	
	\newcommand{\ut}{\ti{u}}
	\newcommand{\ue}{u^{\e}}
	\newcommand{\us}{u_{\sigma}}
	\newcommand{\usx}{(u_{\sigma})_{x}}
	\newcommand{\usxx}{(u_{\sigma})_{xx}}
	\newcommand{\uast}{u_{\sigma}^{\ast}}
	\newcommand{\uen}{u^{\e}_{0}}
	\newcommand{\ua}{u^{\alpha}}
	\newcommand{\uam}{u^{\alpha-1}}
	\newcommand{\uap}{u^{\alpha+1}}
	\newcommand{\ux}{u_{x}}
	\newcommand{\uxx}{u_{xx}}
	\newcommand{\uxxx}{u_{xxx}}
	\newcommand{\uex}{\ue_{x}}
	\newcommand{\uexx}{\ue_{xx}}
	\newcommand{\uexxx}{\ue_{xxx}}
	\newcommand{\un}{u_{0}}
	\newcommand{\uea}{(\ue)^{\alpha}}
	\newcommand{\ueam}{(\ue)^{\alpha-1}}
	\newcommand{\ueap}{(\ue)^{\alpha+1}}
	\newcommand{\lal}{\lambda_{k}}
	\newcommand{\lan}{\lambda_{0}}
	\newcommand{\laml}{\lambda_{-k}}
	\newcommand{\ssl}{\sigma_{s}^{k}}
	\newcommand{\Cstr}{C_{\text{Strat}}}
	\newcommand{\xiast}{\xi^{\ast}}
	\newcommand{\fast}{f^{\ast}}
	
	\newcommand{\uet}{\ti{u}^{\e}}
	\newcommand{\ps}{p_{\sigma}}
	\newcommand{\psx}{(\ps)_{x}}
	\newcommand{\uent}{\ti{u}^{\e}_{0}}
	\newcommand{\uat}{\ti{u}^{\alpha}}
	\newcommand{\uamt}{\ti{u}^{\alpha-1}}
	\newcommand{\uapt}{\ti{u}^{\alpha+1}}
	\newcommand{\uxt}{\ti{u}_{x}}
	\newcommand{\uxxt}{\ti{u}_{xx}}
	\newcommand{\uext}{\uet_{x}}
	\newcommand{\uexxt}{\uet_{xx}}
	\newcommand{\unt}{\ti{u}_{0}}
	\newcommand{\ueat}{(\uet)^{\alpha}}
	\newcommand{\ueamt}{(\uet)^{\alpha-1}}
	\newcommand{\ueapt}{(\uet)^{\alpha+1}}
	\newcommand{\Gas}{G_{\alpha }}
	\newcommand{\Gmcas}{\mathcal{G}_{\alpha }}
	\newcommand{\wt}{\ti{w}}
	\newcommand{\vt}{\ti{v}}
	\newcommand{\zt}{\ti{z}}
	\newcommand{\ve}{v^{\e}}
	\newcommand{\ze}{z^{\e}}
	\newcommand{\vet}{\ti{v}^{\e}}
	\newcommand{\wet}{\ti{w}^{\e}}
	\newcommand{\zet}{\ti{z}^{\e}}
	\newcommand{\vext}{\ti{v}^{\e}_{x}}
	\newcommand{\zext}{\ti{z}^{\e}_{x}}
	\newcommand{\vexxt}{\ti{v}^{\e}_{xx}}
	\newcommand{\zetxx}{\ti{z}^{\e}_{xx}}
	
	\newcommand{\epsw}[1]{\e\, \cW(#1)}
	\newcommand{\epsws}[1]{\e\, \cW'(#1)}
	\newcommand{\epswss}[1]{\e\, \cW''(#1)}
	\newcommand{\gl}{g_{k}}
	\newcommand{\gml}{g_{-k}}
	\newcommand{\glx}{(g_{k})_{x}}
	\newcommand{\glxx}{(g_{k})_{xx}}
	
	\newcommand{\bl}{\beta_{k}}
	\newcommand{\blt}{\ti{\beta}_{k}}
	\newcommand{\ble}{\beta^{\e}_{k}}
	\newcommand{\blet}{\ti{\beta}^{\e}_{k}}
	\newcommand{\W}{F}
	\newcommand{\We}{W^{\e}}
	\newcommand{\Wt}{\ti{W}}
	\newcommand{\Wet}{\ti{W}^{\e}}
	\newcommand{\Wpq}{\mathcal{W}_{pq\text{Strat}}}
	\newcommand{\WStrat}{F_{\text{Strat}}}
	\newcommand{\Tmax}{T}
	\newcommand{\Ot}{\cO_{T}}
	\newcommand{\Otb}{\bar{\cO}_{T}}
	
	\newcommand{\Norm}[1]{\lVert \,  #1 \, \rVert }
	\newcommand{\Normsq}[1]{\lVert \, #1 \, \rVert^{2} }
	\newcommand{\skp}[2]{\left( #1  ,  #2  \right)}
	\newcommand{\iprod}[2]{\langle \,  #1 \, , \,#2 \, \rangle}
	
	\newcommand{\abl}[1]{\frac{\partial}{\partial #1}}
	\newcommand{\ablsq}[1]{\frac{\partial^{2}}{\partial #1^{2}}}
	\newcommand{\ablcu}[1]{\frac{\partial^{3}}{\partial #1^{3}}}
	
	\newcommand{\stochbas}{(\Omega, \mathcal{F}, (\mathcal{F}_{t})_{t \ge 0},\mathbb{P})}
	\newcommand{\stochbasti}{(\ti{\Omega}, \ti{\mathcal{F}}, (\ti{\mathcal{F}_{t}})_{t \ge 0},\ti{\mathbb{P}})}
	\newcommand{\qv}[1]{\langle #1 \rangle}
	\newcommand{\crossv}[2]{\langle\langle #1,#2 \rangle\rangle}
	
	%\theoremstyle{break}
	%\theoremseparator{}
	\newtheorem{sat1}{Theorem}[section]
	\newtheorem{def1}[sat1]{Definition}
	\newtheorem{the1}[sat1]{Theorem}
	\newtheorem{lem1}[sat1]{Lemma}
	\newtheorem{kor1}[sat1]{Corollary}
	\newtheorem{bsp1}[sat1]{Example}
	\newtheorem{bem1}[sat1]{Remark}
	\newtheorem{prop1}[sat1]{Proposition}
	
	\newenvironment{myproof}[1]%
	{\paragraph{\textbf{Proof{#1}:}}}
	%{\par\noindent{\textit{ Proof{#1}:}}\hspace{0.5em}}%
	{\\ \hspace*{\fill}$\square$\newline} 
	
	%Aus FischerGruen
	\def\Xint#1{\mathchoice
		{\XXint\displaystyle\textstyle{#1}}%
		{\XXint\textstyle\scriptstyle{#1}}%
		{\XXint\scriptstyle\scriptscriptstyle{#1}}%
		{\XXint\scriptscriptstyle\scriptscriptstyle{#1}}%
		\!\int}
	\def\XXint#1#2#3{{\setbox0=\hbox{$#1{#2#3}{\int}$ }
			\vcenter{\hbox{$#2#3$ }}\kern-.6\wd0}}
	\def\ddashint{\Xint=}
	\def\dashint{\Xint-}
	
	\newcommand{\expect}{{\mathbb E}}
	\newcommand{\fbeta}{{F_\beta}}
	\newcommand{\vdm}{{|v_\delta|^{m-1}}}
	\newcommand{\meanball}{{\dashint_{B_\delta(x_0)}}}
	\newcommand{\Rball}{{B_R(x_0)}}
	\newcommand{\rball}{{B_r(x_0)}}
	\newcommand{\ort}{{\mathcal O}}
	\newcommand{\intort}{{\int_\ort}}
	\newcommand{\filt}{{\mathcal F_t}}
	\newcommand{\fils}{{\mathcal F_s}}
	\newcommand{\rhod}{{\rho_\delta\ast}}
	\newcommand{\fdu}{{F_\delta[u]}}
	\newcommand{\probab}{{\mathbb P}}
	\newcommand{\tprobab}{{\tilde\probab}}
	\newcommand{\integrand}{{|u|^{2m}(x,t)dxdt}}
	\newcommand{\intte}{{\int_0^{T_E}}}
	
	\newcommand{\Ih}{\mathcal{I}_h}
	\newcommand{\domain}{\mathcal{O}}
	\newcommand{\Lh}{{L_h}}
	%\tableofcontents
	
	\title[Zero-contact angle solutions to stochastic thin-film equations]{Zero-contact angle solutions to stochastic thin-film equations}

	\author[G\"unther Gr\"un and Lorenz Klein]{G\"unther Gr\"un and Lorenz Klein}
	\date{\today}
	\keywords{Stochastic partial differential equation, Degenerate parabolic equation, Fourth-order equation, Thin fluid film; Thermal fluctuations; Stratonovich noise; Martingale solution, Free-boundary problem, Compactly supported initial data}
	\makeatletter
	\@namedef{subjclassname@2020}{%
		\textup{2020} Mathematics Subject Classification}
	\makeatother 
	\subjclass[2020]{60H15, 76A20, 76D08, 35K35, 35K65, 35Q35, 35R35, 35R37} 
% 	60H15   	Stochastic partial differential equations (aspects of stochastic analysis)
%   35R35   	Free boundary problems for PDEs
%   35R37       Moving boundary problems
%   65M12   	Stability and convergence of numerical methods for initial value and initial-boundary value problems involving PDEs
%   65M60   	Finite element, Rayleigh-Ritz and Galerkin methods for initial value and initial-boundary value problems involving PDEs
%   76A20       Thin Film Flow
%   76D08   	Lubrication theory
%   35K35   	Initial value problem for higher-order parabolic equations
%   35K65       Degenerate parabolic pde
%   35Q35   	PDEs in connection with fluid mechanics
%
%	
\begin{abstract}
We establish existence of nonnegative martingale solutions to stochastic thin-film equations with compactly supported initial data under Stratonovich noise. Based on so called $\alpha$-entropy estimates, we show that almost surely these solutions are classically differentiable in space almost everywhere in time and that their  derivative attains the value zero at the boundary of the solution's support. I.e., from a physics perspective, they exhibit a zero-contact angle at the three-phase contact line between liquid,  solid, and ambient fluid. These $\alpha$-entropy estimates are first derived for almost surely strictly positive solutions to a family of stochastic thin-film equations augmented by second-order linear diffusion terms. Using It\^o's formula together with stopping time arguments, the Jakubowski/Skorokhod calculus, and martingale identification techniques, the passage to the limit of vanishing regularization terms gives the desired existence result.
\end{abstract}	
	\maketitle

	\parsep 12pt
	\parskip 5pt
	\parindent 0pt
	
\section{Introduction}
	In this paper, we are concerned with existence results of martingale solutions to stochastic  thin-film equations of the generic form
	\begin{align}\label{eq:STFES}
		\dif u =  - (m(u)u_{xxx})_{x} \dt 
		+ (\sqrt{m(u)} \circ \dif W)_{x}  
	\end{align}
	subject to periodic boundary conditions.
        The deterministic version of \eqref{eq:STFES} models the solely surface-tension driven evolution of the height $u$ of a thin viscous liquid film -- the noise term is to capture effects of thermal fluctuations.
        
	Gess and Gnann have been the first to consider stochastic thin-film equations with Stratonovich noise.
	In \cite{GessGnann2020}, they proved the global-in-time existence of nonnegative martingale solutions for the choice $m(u)=u^{2}$.
	To establish this result they took advantage of the regularizing effect of Stratonovich noise compared to It\^o noise. 
	In fact, for the It\^o version of \eqref{eq:STFES}, i.e.
	\begin{align}
		\dif u =  - (m(u)u_{xxx})_{x} \dt 
		+ (\sqrt{m(u)} \dif W)_{x}  \, , 
	\end{align}
	no integral estimates are known.
	In contrast, the Stratonovich version \eqref{eq:STFES} permits to derive stochastic versions of the energy estimate
	\begin{align}\label{est:energy}
		\frac{1}{2}	\intsp{\abs{ \ux(\cdot,t)}^{2}} + \intsptit{u^{2} u_{xxx}^{2}} \le \frac{1}{2} \intsp{\abs{(\un)_{x}}^{2}}
	\end{align} 
	and the so called entropy estimate
	\begin{align}\label{est:entest}
		\intsp{G(u(\cdot,t))} + \intsptit{u_{xx}^{2}} \le \intsp{G(\un)} \, ,
	\end{align}
	where $G(\cdot)$ is a second primitive of the reciprocal mobility $m^{-1}(s)$.
	Already from the deterministic setting, it is well-known that weak solutions to the free-boundary problem associated with  the thin-film equation are not unique in general, unless additional conditions are imposed at the free boundary, i.e. the boundary of $\operatorname{supp} [u(\cdot, t) > 0]$.
	In a series of papers 
	\cite{GKO2008,GK2010,GGKO2014,Gnann2015,Gnann2016,GIM2019}
	short-time uniqueness results were established for classical solutions of thin-film equations exhibiting a zero-contact angle at the free boundary.
        \newline\newline
	In this spirit, it is the aim of the present paper to construct nonnegative martingale solutions $\ut$ to equation \eqref{eq:STFES} which are $\ti{\mathbb{P}}$-almost surely and almost everywhere in time continuously differentiable in space. Hence, those spatial derivatives of $\ut$ attain the value zero  in roots of $\ut$. 
	The vanishing of these derivatives comes as the consequence of additional regularity results. 
	While the solutions constructed by Gess and Gnann in \cite{GessGnann2020} for compactly supported initial data do not have the regularity stipulated by the entropy estimate \eqref{est:entest} (note that $\int_{\cO}G(\un) \dx = + \infty$ in this case) and therefore do not necessarily exhibit zero contact angles, the solutions presented here are more regular. 
	In fact, they satisfy a stochastic version of a variant of \eqref{est:entest}, the so-called $\alpha$-entropy estimate. 
	This $\alpha$-entropy estimate provides $H^{2}$-regularity of appropriate powers of the solutions $\ut$ without requiring initial data to be zero only on sets of Lebesgue measure zero.
        For an overview on $\alpha$-entropy estimates and other integral estimates for the thin-film equation in the deterministic setting, we refer to \cite{Beretta1995,BP1996, DalPassoGarckeGruen, GruenDroplet} and the references therein.
	
	At this point it is worth mentioning that in the analysis of the qualitative behaviour of deterministic thin-film equations, weighted versions of $\alpha$-entropy estimates become important. They have been used, e.g.,  to obtain optimal results on the propagation of the free boundary of solutions or on the regularity at the free boundary.
	
	For an overview of corresponding results, we cite 
	\cite{GruenDropletPropagRate,BertschDalPassoGarckeGruen, FischerOptRate} for finite speed of propagation and \cite{DalPassoGiacoGruenWTP, GiacomelliGruenWTP, GruenUpperBounds,FischerUpperBoundWTP} for the occurrence and scaling of waiting time phenomena. 
	In the stochastic setting, the techniques of \cite{GiacomelliGruenWTP} have been generalized by \cite{FischerGruenFSP,GrillmeierDiss} to provide sufficient criteria for the occurrence of waiting time phenomena and for qualitative results on finite speed of propagation for stochastic $p$-Laplace and stochastic porous-media equations. For finite speed of propagation for the latter equations, we also mention \cite{BarbuRoeckner2012, GessFSP2013} which use different techniques.
	
	Before giving the outline of the present paper, we report on variants of \eqref{eq:STFES} which are meaningful for physical and/or for analytical reasons as they may set up auxiliary problems to construct the more regular solutions to be considered in this paper.
	First, we mention \eqref{eq:STFES} with the generic mobility $m(u)= u^{n}$ where $n>0$.
	The exponent $n$ depends on the flow boundary conditions at the liquid-solid interface -- a no-slip boundary conditions entails $n=3$.
	Recently, Dareiotis, Gess, Gnann, and the first author of this paper established the existence of martingale solutions \cite{DGGG} for \eqref{eq:STFES} with $m(u)=u^{n}$ in the parameter regime $n \in [8/3,4)$ which covers in particular the no-slip case.	
	Note that Davidovitch et al. \cite{DavidovitchStone2005} who derived \eqref{eq:STFES} with It\^o- instead of Stratonovich noise via the dissipation-fluctuation theorem conjectured that noise enhances spreading, changing in particular characteristic spreading laws on intermediate time-scales in expectation.    
	Parallel in time, Gr\"un, Mecke, and Rauscher \cite{MeckeRauscher} studied the influence of thermal fluctuations on the dewetting of unstable liquid films.
	Based on lubrication approximation and Fokker-Planck-type arguments, they came up with an equation of the generic form
	\begin{align}\label{eq:STFE}
		\dif u =  - (m(u) (\uxx - \W'(u))_{x})_{x} \dt 
		+ (\sqrt{m(u)} \dif W)_{x} \, ,  
	\end{align}
	where the effective interface potential $\W(u)$ models van der Waals-interactions -- a typical example is the potential $\W(u) := \alpha u^{-p}- \beta u^{-q}$ with $p>q>0, \alpha >0$, and $\beta \ge 0$. 
	For the case $m(u)=u^{2}$, the existence of a.s. positive martingale solutions has been established in \cite{FischerGruen} -- the technically much more involved case of two space dimensions has been studied in \cite{MetzgerGruen2021} -- for a very recent result in the spirit of \cite{GessGnann2020} which also provides $\alpha$-entropy estimates, see \cite{Sauerbrey2021}.
	
	The outline of our paper is as follows.
	In contrast to the Trotter-Kato scheme, where the stochastic and the deterministic parts of the equation are split and which was used in \cite{GessGnann2020}, we will follow an approximation ansatz based on positive solutions to 
	\begin{align}\label{eq:STFEFGapprox}
		\dif \ue = - ( m(\ue) ( \ue_{xx} - \e\W'(\ue))_{x})_{x} \dt + (\sqrt{m(\ue)}\circ \dif W)_{x} \, ,
	\end{align}
	where $\e \in (0,1)$ and $F(u) := u^{-p}$, $p>2$. 
%	The benefit of our approach is that by means of Fatou's lemma and arguments based on lower semi-continuity the regularity of the approximate solutions can be recovered in the limit.
	We take advantage of the fact that under natural assumptions on the coloured noise $W(x,t) = \sum_{k \in \Z} \lal \gl(x) \bl(t)$, such that
	\begin{itemize}
		\item
		$
		 \gl(x) = 
		\begin{cases}
		\sqrt{\frac{2}{L}} \sin(\frac{2 \pi k x }{L})\quad &k >0, 
		\\
		\frac{1}{\sqrt{L}}\quad &k=0 ,
		\\
		\sqrt{\frac{2}{L}} \cos(\frac{2 \pi k x }{L})\quad &k<0,
		\end{cases}
		$
		\item
		$\bl(\cdot)$, $k \in \Z$, are i.i.d. Brownian motions on $\R$,
		\item 
		  $\lal \in \R_{0}^{+}$, $k \in \Z$, appropriate damping
                  parameters,  
	\end{itemize} 
	(see Section \ref{sec:Preliminaries} for the precise assumptions), existence of a.s. positive martingale solutions to \eqref{eq:STFEFGapprox} comes as a consequence of the existence result in \cite{FischerGruen}.
	
	After formulating precise assumptions in Section~\ref{sec:Preliminaries}, especially on  initial data and the noise, 
	we present our main results in Section~\ref{sec:mainresult}.
	The existence of  solutions to \eqref{eq:STFEFGapprox} is the topic of
        Section~\ref{sec:ApproxSolutions}. 
	Our strategy  is the following. 
	Since we are dealing with Stratonovich noise in \eqref{eq:STFEFGapprox}, we may rewrite it in It\^o form with the corresponding correction term added. 
	The equation obtained is then expressed in the form \eqref{eq:STFE},
just including the Stratonovich correction term in the potential $F$, cf. 
 \eqref{eq:stratito1},  \eqref{eq:stratito2} and \eqref{eq:FischerGruenModified}. That way,  the equation satisfies -- with a grain of salt -- the assumptions of the existence result, Theorem~3.2,   in \cite{FischerGruen}.
		Since solutions in \cite{FischerGruen}  were constructed under the assumption of positive initial data, we shift nonnegative initial data with potentially compact support by a suitable power of $\e$ such that we recover the nonnegative initial data in the limit, cf. (H2$\e$) in Section \ref{sec:ApproxSolutions}.
        This is sufficient to establish in Theorem~\ref{theo:existence-stfe-stratonovich} existence results for a family of approximate $\ti\probab$-almost surely strictly positive  solutions $u^\e.$  
	
	The key result for the passage to the limit $\e\to 0$ is a combined $\alpha$-entropy-energy estimate in the spirit of the classical $\alpha$-entropy estimates in \cite{Beretta1995} (see also \cite{BP1996}) and the energy estimates in \cite{Bernis1990} both translated to the stochastic setting.
	The derivation of this estimate is the content of the fifth section. 
 	We first introduce suitable stopping times and cut-off versions of our approximate solutions, cf. \eqref{def:stoppingTimes} and \eqref{def:ModifiedSolutions}, which allow  to derive a first version of an $\alpha$-entropy-energy estimate  in Theorem \ref{lem:alpha-Entropy-Energy}.
	
	It\^o's formula, which is the main tool for the proof, is applied to the energy $\int_\ort u_x^2dx$ and the $\alpha$-entropy $\int_\ort \tfrac{1}{\alpha(\alpha+1)} u^{\alpha+1}-\tfrac{1}{\alpha} u +\tfrac{1}{\alpha+1}dx$, see Appendix \ref{appendix:Itoformula} for the rigorous justification.
	Here, the advantages of the usage of the Stratonovich integral become apparent again. Critical terms occurring in  It\^o's formula are controlled by
         the Stratonovich correction term -- this way guaranteeing the estimate to be $\e$-independent.

	The passage to the limit $\e\to 0$  is discussed in Section \ref{sec:Jakubowski}. 
	We use the aforementioned $\alpha$-energy-entropy estimate to apply Jakubowski's theorem, cf. \cite{Jakubowski1994}.
	Based on this, we follow standard arguments encountered in the analysis of PDE's for the convergence of the deterministic terms and make use of the ideas introduced in \cite{Brezniak2007, HofmanovaWeak2016} to identify the stochastic integral in Lemma \ref{sat:IdentM0}. 
	The effective interface potential vanishes in the limit, cf. Lemma \ref{lem:PotentialTermVanishes}. 
	
	There are three appendices. In Appendix \ref{StratCorrection} the equivalence of the different formulations of \eqref{eq:STFEFGapprox} is made explicit.
	More details on the application of It\^o's formula are provided in appendix \ref{appendix:Itoformula} and \ref{predictability}.  	
	
	\textbf{Notation:} 
	Besides the standard notation of pde theory and stochastic analysis, we use the following. 
	By $C$ we  denote a generic constant. 
	Throughout the paper we will use $\e$ as an approximation parameter, subsequences will not be renamed, if it causes no confusion.     
	We consider the spatial domain $\cO = [0,L]$ and define $\cO_{T} := (0,L) \times (0,T)$ for numbers $L,T > 0$. 
	For a function $f$ on $\Ot$, $[f>0]$ is the set $\{(x,t) \in \Ot; \, \, f(x,t) >0	\}$.
	Subspaces consisting of periodic functions (w.r.t. space) are marked by the subscript `per' on the corresponding function space. 
	For $\gamma,\sigma \in (0,1] $ we denote by $C^{\gamma,\sigma} 
	(\Otb)$ the space of Hölder-continuous functions on $[0,L]$ w.r.t. the exponent  $\gamma$ and on $[0,T]$ w.r.t. $\sigma$.
	The minimum of $a$ and $b$ is denoted by $a \wedge b$. 
	We write $\qv{X}$ for the quadratic variation process of a stochastic process $X$ and $\qv{X,Y}$ for the quadratic covariation process of $X$ and another process $Y$. 
	Moreover, for two Hilbert spaces $U$ and $V$, $L_{2}(U,V)$ is the set of Hilbert-Schmidt operators from $U$ to $V$. Note that the dual pairing on a Banach space $X$ is denoted by ${}_{X'}\langle x',x\rangle_{X}$ for $x'\in X'$ and $x\in X$.

\section{Preliminaries}\label{sec:Preliminaries}
	
	Let us fix some basic assumptions. We are dealing with the stochastic thin-film equation with Stratonovich noise
	\begin{align}\label{eq:StfeStrato}
		\dif u =  - \left(m(u)u_{xxx}\right)_{x} \dt 
		+ \left(\sqrt{m(u)} \circ \dif W\right)_{x}  
	\end{align}
	on $\Ot$ subject to periodic boundary conditions and initial data specified below.
	For the noise we consider a Q-Wiener process defined by the operator 
	\begin{align}\label{KovarOp}
		Q\gl = \lal^{2}\gl \quad \forall k \in \Z  .
	\end{align}
	Here, the functions $\gl$ form a basis of $\lpr{2}{\cO}$ consisting of eigenfunctions of the Laplacian on $\cO$ subject to periodic boundary conditions.
	\begin{align}
		\label{def:ONB} \gl(x) = 
		\begin{cases}
			\sqrt{\frac{2}{L}} \sin(\frac{2 \pi k x }{L})\quad &k >0 
			\\
			\frac{1}{\sqrt{L}}\quad &k=0 
			\\
			\sqrt{\frac{2}{L}} \cos(\frac{2 \pi k x }{L})\quad &k<0
		\end{cases}
	\end{align}
	The noise is coloured by the growth condition on the numbers $(\lal)_{k \in \Z}$, cf. (H3) below. 
	We can now give precise assumptions for our main result. 
	\begin{itemize}
		\item[(H1)] 
		The mobility is given by $m(u) = u^{2}$.
		\item[(H2)]
		Let $\Lambda^{0}$ be a probability measure on $H^{1}_{\text{per}}(\cO)$ equipped with the Borel $\sigma$-algebra which is supported on the subset of nonnegative functions such that there is a positive constant $C$ with the property that 
		\begin{align*}
			\operatorname{esssup}_{v \in \operatorname{supp \Lambda^{0}}} 
			\left\{		
			\int_{\cO} \frac{1}{2} \abs{v_{x}}^{2}\dx + \left( \int_{\cO} v \dx \right)   	
			\right\}
			\le C \, .
		\end{align*}
		\item[(H3)]
		Let $\stochbas$ be a stochastic basis with a complete right-continuous filtration such that
		\begin{itemize}
			\item[-] 
			$W$ is a $Q$-Wiener process on $\Omega$ adapted to $(\mathcal{F}_{t})_{t \ge 0}$ which admits a decomposition of the form 
			$W = \sum_{k \in \Z} \lal \gl \bl$ 
			for a sequence of independent standard Brownian motions $\bl$ and nonnegative numbers $(\lal)_{k \in \Z}$ with 
			\begin{align}\label{def:Eigenvalues}
				\laml = \lal 
			\end{align}
			for all  $k \in \N$,
			\item[-]
			the noise is coloured in the sense that 
			\begin{align}\label{Lambda-Wachstumsbedingung}
				\sum_{k \in \Z} k^{4} \lal^{2} < \infty \,, 
			\end{align}	
			\item[-]
			there exists a $\mathcal{F}_{0}$-measurable random variable $\un$ such that  $\Lambda^{0} = \mathbb{P} \circ \un^{-1}$.
		\end{itemize}				
	\end{itemize}
	Based on these hypotheses, we may rewrite \eqref{eq:StfeStrato} in two different ways.
	\begin{equation}\label{eq:stratito1}
		\dif u = -( u^2(u_{xx}-\mathcal S'(u))_x)_x\dif t +(u\dif W)_x
	\end{equation}
	with $\mathcal S(u):= C_{Strat}\left(u-\log u\right)$ and $C_{Strat}:=\tfrac12\left(\tfrac{\lambda_0^2}{L}+\sum_{k=1}^\infty \tfrac{2\lambda_k^2}{L}\right)$.
	Note that this is equivalent to 
	\begin{equation}\label{eq:stratito2}
		\dif u=(-(u^2u_{xxx})_x+C_{Strat} u_{xx}) \dif t +(u \dif W)_x \, .
	\end{equation}
        For a justification, we refer to Appendix A.

\section{Main Results}\label{sec:mainresult}
In this section, we make our results on the existence of zero-contact angle martingale solutions precise.

	\begin{the1}\label{theo:MainResult}
		Let (H1), (H2), and (H3) be satisfied and let $T>0$ be given. 
		Then there exist a stochastic basis $\stochbasti$ and a $\ti{\mathcal{F}}_{t}$-adapted Q-Wiener process $\ti{W} = \sum_{k \in \Z}\lal\gl \blt$, 
		a stochastic process $\ut \in \lpr{2}{\ti{\Omega};L^{2}(0,T;W^{1,3}_{\text{per}}(\cO))} \cap L^{2}(\ti{\Omega}; C^{\ti{\gamma}, \ti{\sigma}}(\Otb))$, $\ti{\gamma}<1/2, \, \ti{\sigma} < 1/8$, and $\unt \in \lpr{2}{\ti{\Omega};H^{1}_{\text{per}}(\cO)}$
		such that the following holds:
		\begin{enumerate}
			\item \label{theo:MainResult1}
			$\ut$ and $\unt$ are $\ti{\mathbb{P}}$-almost surely nonnegative,
			\item\label{theo:MainResult3}
			for $t \in [0,T]$ and all $\phi \in H^{3}_{\text{per}}(\cO)$
			\begin{align}
				\label{eq:MainResultWeakForm}
				\intsp{(\ut(t)-\unt) \phi} 
				&= \int \int_{[\ut > 0]}\ut_{x}^3\phi_{x} \dx \ds \notag
				+ 3 \intsptit{\ut \ut_{x}^{2}\phi_{xx}} 
				\\ \notag
				&+ \intsptit{\ut^{2}\ut_{x}\phi_{xxx}}	
				\\ \notag
				&- \frac{1}{2} \intsptit{\sum_{k\in \Z} \lal^{2}\gl(\gl\ut)_{x}\phi_{x}} 
				\\ 
				&- \sum_{k \in \Z}\int_{0}^{t}\intsp{ \lal\gl\ut\phi_{x}} \dblt
			\end{align}  
			holds $\ti{\mathbb{P}}$-almost surely and we have  
			$
			\Lambda^{0} = \ti{\mathbb{P}} \circ \unt^{-1},
			$
			\item \label{theo:MainResult4}
			$\ut$ satisfies for arbitrary $q \ge 1$ and $\alpha \in (-\frac{1}{3},0)$ the estimate  
			\begin{align}\label{est:MainResult} \notag	 
				\Ew{\sup_{t\in [0,T]} \left (\intsp{\frac{1}{2} \abs{\ut_{x}}^{2}}  \right)^{q}} 
				&+
				\Ew{\left(\intsptiT{((\ut)^{\frac{\alpha + 3}{4}})_{x}^{4}}\right)^{q}}
				\\
				&+ 
				\Ew{\left(\intsptiT{((\ut)^{\frac{\alpha+3}{2}})^{2}_{xx}}\right)^{q}} 
				\le  C(\unt,q,T).
			\end{align}
		\end{enumerate}
	\end{the1}
\begin{kor1}\label{cor:contactangle}
Let $\ut$ be a solution as constructed in Theorem~\ref{theo:MainResult} and $\alpha\in (-1/3,0).$ Then, 	$\ti{\mathbb{P}}$-almost surely, $\ut$ exhibits a zero-contact angle in the following sense: For almost all $t_0\in (0,T]$, the classical derivative $\frac{\partial}{\partial x} \ut(x_0,t_0,\omega)$ exists in points $x_0\in\ort$ such that $\ut(x_0,t_0,\omega)=0$, and it attains the value zero.
\end{kor1}

\section{An Existence Result for Positive Approximate Solutions}
	\label{sec:ApproxSolutions}
        As pointed out in the introduction, our existence result for zero-contact angle solutions relies on new integral estimates which are initially derived for strictly positive approximate solutions. In this section, we provide  results on existence  and on refined regularity and positivity properties of such solutions. More precisely, we study stochastic thin-film equations
	\begin{equation}\label{eq:FischerGruenModified} 
		\dif \ue =  -((\ue)^2(\ue_{xx} - \Pi_\varepsilon'(\ue))_{x})_{x} \dt  
		+ \sum_{k \in \Z} \lal (\ue\gl)_{x}  \dble 	,
	\end{equation} 
	$\varepsilon\in(0,1)$, subject to periodic boundary conditions on $\cO$. This class of equations differs from equations \eqref{eq:stratito1} and \eqref{eq:stratito2} by the choice of the potential $\mathcal S$ (or $\Pi_\varepsilon$, respectively) and by the assumptions on the positivity properties of initial data.
	Here, $\Pi_\varepsilon(u)$ is given by
	\begin{equation}
		\label{eq:defPi}
		\Pi_\varepsilon (u):=\begin{cases}\varepsilon u^{-p} + \mathcal S(u) & \text{ if } u>0\\
			+\infty & \text{ if } u\leq 0
		\end{cases}	
	\end{equation}
	with a positive number $p>2.$
	Note that the potential $\Pi_\varepsilon$ satisfies for every $\varepsilon>0$ the hypothesis	
	\begin{itemize}
		\item[(H4$\e$)]\label{H4new}
		For $\e>0$, the effective interface potential $\Pi_{\e}$ has continuous second-order derivatives on $\R^{+}$ and satisfies for some $p > 2$ and $u > 0$
		\begin{align*}
			\notag
			c_{1} u^{-p-2} - c_{2} &\le \Pi''_{\e}(u) \le C_{1} (u^{-p-2}+1) \, ,
			\\
			&\Pi_{\e}(u) \ge Cu^{-p} \,  ,
		\end{align*}
		where $c_{1}, C_{1}$, and $C$ are positive constants depending on $\e>0.$
	\end{itemize}
	In addition, we require initial data to be bounded from below by $\e^\theta$ with $\theta\in (0, 1/p)$. Combined with Hypothesis (H2), this leads to the modification 
	\begin{itemize}
		\item[(H2$\e$)]
		Let $\Lambda_0$ be a probability measure on $H^1_{\text{per}}(\ort)$ satisfying Hypothesis (H2). Then, $\Lambda^{\e}$ is the probability measure on $H^1_{\text{per}}(\ort)$ defined by 
		$\Lambda^{\e} = \Lambda^{0} \circ S_{\e}^{-1}$, where \newline $S_{\e}: H^{1}_{\text{per}}(\cO) \rightarrow H^{1}_{\text{per}}(\cO)$ is 
		for $\theta \in (0,1/p)$ 
		given by $S_{\e}(u) = u + \e^{\theta}$. 	
	\end{itemize}
	
	\begin{def1}\label{def:weak-martingale-solution}
		Let $\Lambda^{\e}$ be a probability measure on $H^{1}_{\text{per}}({\cO})$  satisfying $\text{(H2$\e$)}$. We call a triple $((\Omega^{\e}, \mathcal{F}^{\e}, (\mathcal{F}_{t}^{\e})_{t\ge0},\mathbb{P}^{\e}),\ue,W^{\e})$ a  martingale solution to \eqref{eq:FischerGruenModified} with initial data $\Lambda^{\e}$ on the time interval $[0,T]$ provided 
		\begin{itemize}
			\item[i)]
			$(\Omega^{\e}, \mathcal{F}^{\e}, (\mathcal{F}_{t}^{\e})_{t\ge0},\mathbb{P}^{\e})$
			is a stochastic basis  with a complete right-continuous filtration,
			\item[ii)]
			$W^{\e}$ satisfies (H3) with respect to $(\Omega^{\e}, \mathcal{F}^{\e}, (\mathcal{F}_{t}^{\e})_{t\ge0},\mathbb{P}^{\e})$,	
			\item[iii)]
			$\ue \in L^{2}(\Omega^{\e}; L^{2}(0,T;H^{3}_{\text{per}}(\cO))) \cap L^{2}(\Omega^{\e}; C^{\gamma ,\sigma}(\Otb))$ with $\gamma< 1/2,\,\sigma < 1/8 $ is $\mathbb{P}^{\e}$-almost surely positive and adapted to $(\mathcal{F}_{t}^{\e})_{t\ge0}$,		
			\item[iv)]
			there is  a $\mathcal{F}_{0}^{\e}$-measurable $H^{1}_{\text{per}}(\cO)$-valued random variable $\ue_{0}$ with $ \Lambda^{\e} = \mathbb{P}^{\e} \circ (\uen)^{-1}$
			and the equation
			\begin{align}\label{eq:SolutionApproxProblems1}
				\intsp{(\ue(t)-(u_{0}^{\e}))\phi} \notag
				&= \intsptit{m(\ue)(\ue_{xx}- \Pi_{\e}'(\ue))_{x} \phi_{x} } 
				\\ 
				&- \sum_{k \in \Z} \int_{0}^{t} \int_{\cO} \lal \gl\sqrt{m(\ue)}\phi_{x} \dx \dble
			\end{align}
			holds $\mathbb{P}^{\e}$-almost surely for all $t\in [0,T]$ and all $\phi \in H^{1}_{\text{per}}(\cO) $.
		\end{itemize}
	\end{def1}
	We  have the following theorem.
	\begin{the1}\label{theo:existence-stfe-stratonovich}
		Let the assumptions (H1), (H2$\e$), (H3), and (H4$\e$) be satisfied  and let $T >0$ be given.
		Then, for every $\e \in (0,1)$ there exists a  martingale solution \\ $((\Omega^{\e}, \mathcal{F}^{\e}, (\mathcal{F}_{t}^{\e})_{t\ge0},\mathbb{P}^{\e}),\ue,W^{\e})$ to \eqref{eq:FischerGruenModified} in the sense of Definition \ref{def:weak-martingale-solution}, satisfying the additional bound
		\footnote{Here and in what follows expectations of terms that depend on $\e$ are computed w.r.t. the probability measure $\mathbb{P}^{\e}$.} 
		\begin{align}\notag
			\label{est:EnEntropFischeGruen}
			&\Ew{\sup_{t\in [0,\Tmax]} 
				\left( 
				\int_{\cO}\frac{1}{2} \abs{\uex}^{2} + \Pi_{\e}(\ue)	\dx 
				\right)^{q} } 
			+
			\Ew{ 
				\int_{0}^{T} \int_{\cO} (\ue)^{2} \abs{(\uexx -  \Pi'_{\e}(\ue))_{x}}^{2} \dx \dt 
				} 
			\\
			&\le
			C(\e,\un, \Tmax, q) \, ,
		\end{align}  
		where $q \ge 1$ can be chosen arbitrarily.
	\end{the1}
	\proof
	Observe that Hypothesis (H4$\e$) differs only by the positive additive term $C_1$ in the upper bound on $\Pi''_{\e}$ from Hypothesis (H2) in \cite{FischerGruen}. It is straightforward to show that all the results in \cite{FischerGruen} are still true if Hypothesis (H2) of \cite{FischerGruen} is replaced by Hypothesis (H4$\e$) of this paper. Noting that Hypothesis (H3) of this paper is stronger than Hypothesis (H4) of \cite{FischerGruen}, we may establish the result as a consequence of Theorem~3.2 in \cite{FischerGruen}.
	\proofend   
	For brevity we will also call the process $\ue$ in Theorem \ref{theo:existence-stfe-stratonovich} a solution of \eqref{eq:FischerGruenModified}.
	We now give two auxiliary results we will need later on. The first one is a  
	continuous version of Lemma 4.1 in \cite{FischerGruen}. It provides lower bounds on the solutions $\ue$ constructed in Theorem~\ref{theo:existence-stfe-stratonovich}.
	
	\begin{lem1}\label{lem:Minimum-ue}
		Let $\e \in (0,1)$ and consider for $u\in H^1(\ort,\R^+)$ the functional
		\begin{align}
			H_ {\e}[u] := \frac{1}{2}\int_{\cO}  \betr{\ux}^{2} + \e (u)^{-p} \dx \, .
		\end{align}
		Then there is a positive constant $C_p$ independent of $\e$ such that  
		\begin{align}
			\underset{x \in \cO}{\sup}  (u)^{-1} \le \left(\fint_{\cO} u \dx \right)^{-1} + C_{p} \e^{\frac{1}{2-p}} H_{\e}[u]^{\frac{2}{p-2}}.
		\end{align}
		If in addition $H_{\e}[u] \le \sigma^{-1}$ for $ \sigma \in (0,1)$, there is a positive constant $\bar C_p$ independent of $\e$ such that 
		\begin{align}
			\min_{x \in \cO} u \ge \bar{C}_{p} \, \e^{\frac{1}{p-2}} \sigma^{\frac{2}{p-2}} \, .
		\end{align}
	\end{lem1}	
	\proof
	We have 
	\begin{align*}
		\sqrt{\e} \int_{\cO} \abs{\left( (u)^{-\frac{p}{2}+1} \right)_{x}} \dx
		&=
		C(p)\sqrt{\e} \int_{\cO} (u)^{-\frac{p}{2}} \abs{u_x} \dx 
		\\
		&\le
		C(p) 
		\left(
		\e \int_{\cO} (u)^{-p} \dx 
		\right)^{1/2} 
		\left( 
		\int_{\cO} 	\abs{u_x}^{2} \dx 
		\right)^{1/2} 
		\\
		&\le
		C(p)  H_{\e}[u]\, .
	\end{align*}
	This implies
	\begin{align*}
		(\min_{x \in \cO}u)^{-1} 
		= \underset{x \in \cO}{\sup} ((u)^{-1}) 
		= 
		\left(
		\underset{x \in \cO}{\sup} (u)^{-\frac{p-2}{2}}
		\right)^{\frac{2}{p-2}}
		&\le
		\left[ 
		\underset{x \in \cO}{\inf} (u)^{-\frac{p-2}{2}} 
		+
		\e^{-\frac{1}{2}}C(p) 
		\left(
		H_{\e}(u) 
		\right)	
		\right]^{\frac{2}{p-2}} 
		\\
		&\le
		\underset{x \in \cO}{\inf} (u)^{-1} + \e^\frac{1}{2-p}C(p) 
		\left(
		H_{\e}(u) 
		\right)^{\frac{2}{p-2}}
		\\
		&\le
		\left(
		\fint_{\cO} u \dx
		\right)^{-1} 
		+ 
		C_{p} 	\e^{\frac{1}{2-p}} 
		H_{\e}(u)^{\frac{2}{p-2}} \, . 
	\end{align*}
	The assumption  $H_{\e}[u ]\le \sigma^{-1}$ for $\sigma \in (0,1]$ yields
	\begin{align*}
		\min_{x \in \cO} u
		&\ge \frac{1}{
			\left(
			\fint_{\cO} u \dx	
			\right)^{-1}		
			+
			C_{p} \e^{\frac{1}{2-p}} \sigma^{\frac{2}{2-p}}
		}
		\\
		&=
		\frac{\e^{\frac{1}{p-2}} \sigma^{\frac{2}{p-2}}} {
			C_{p}
			+ \left(\e^{\frac{1}{p-2}} \sigma^{\frac{2}{p-2}}\right) \left(\fint_{\cO} u \dx\right)^{-1}
		}
		\\
		&\ge
		\bar{C}_{p} \, \e^{\frac{1}{p-2}} \sigma^{\frac{2}{p-2}} \, .
	\end{align*}
	\proofend	
	In the next lemma, we collect a number of integral estimates which come as natural consequences of the energy-entropy estimate Lemma~4.4 in \cite{FischerGruen}. In the analysis presented in that paper, they were not needed and hence they were not made explicit there.		
	\begin{lem1}\label{lem:UniformEstimate-u^2p}
		Let $\ue$ be a solution as constructed in Theorem~\ref{theo:existence-stfe-stratonovich}. Then there is a positive constant $C=C(\e)$ such that
		\begin{align}\label{lem:UniformEstimate-u^2ph1}
			\Ew{\int_{0}^{T} \int_{\cO} ((\ue)^{2} p^{\e}_{x})^{2} \dx \dt}
			\le
			C
			\,,
		\end{align}
		\begin{align}\label{lem:UniformEstimate-u^2ph2}
			\Ew{\int_{0}^{T} \int_{\cO} \abs{\uexx}^{2} \dx \dt}
			\le
			C
			\,,
		\end{align}
		\begin{align}\label{lem:UniformEstimate-u^2ph3}
			\Ew{\int_{0}^{T} \int_{\cO} (\ue)^{-p-2} (\uex)^{2} \dx \dt}
			\le
			C
			\,.
		\end{align}
	\end{lem1} 
	\proof
	We recall that $\ue$ has been constructed as a limit of solutions $(u^{\e})^{h}$ to finite dimensional auxiliary problems, cf. Lemma 4.2 in \cite{FischerGruen}.  
	The convergence of these solutions is based on the  a-priori estimate in Lemma~4.4 of that paper which will be for us the starting point to derive the estimates above. Adopting    
	the notation of \cite{FischerGruen}, and for the ease of presentation  omitting the index $\e$, we recall the $h$ independent estimate
	\begin{align}\label{lem:UniformEstimate-u^2ph4} \notag
		\Ew{\int_{0}^{T}\int_{\cO} (M_{h}(u^{h}) p^{h}_{x})^{2} \dx \ds  }
		&\le
		C \Ew{\int_{0}^{T} \underset{x\in \cO}{\sup}((u^{h})^{2}) \int_{\cO} M_{h}(u^{h})(p^{h}_{x})^{2} \dx \ds	}
		\\ 
		&\le
		C \Ew{\int_{0}^{T} R(s) \int_{\cO} M_{h}(u^{h})(p^{h}_{x})^{2} \dx \ds	}
		\le
		C(\e, \un,T) \, ,
	\end{align}
	where $R(s) $ denoted in \cite{FischerGruen}  a weighted sum of energy and entropy (cf. (4.14) in \cite{FischerGruen}). In particular, we have the upper bound $\int_\ort (u^h_x)^2(\cdot,s)dx\leq R(s)$ for all $s\in[0,T]$.
	Mimicking the arguments from  Lemma 5.2 in \cite{FischerGruen},  this bound shows that the sequence $(M_{h}(u^{h})p_{x}^{h})_{h \in (0,1)}$ is tight on the path-space $\lpr{2}{\cO \times [0,T]}$. Consequently, by Jakubowski's theorem, cf. \cite{Jakubowski1994}, we find for each $h$ random variables $\ti{g}^{h}: \Omega \rightarrow \lpr{2}{\cO \times [0,T]}$ as well as $\ti{g}: \Omega \rightarrow \lpr{2}{\cO \times [0,T]}$ such that
	\begin{align}\label{lem:UniformEstmate-u^2ph2}
		\ti{g}^{h} \overset{\mathcal{D}} = M_{h}(u^{h})p^{h}_{x}
	\end{align}  under an appropriate probability measure $\ti{\mathbb{P}}$. Moreover, 
	\begin{align*}
		\ti{g}^{h} \rightharpoonup \ti{g}
	\end{align*}
	in $\lpr{2}{\cO \times [0,T]}$ pointwise $\ti{\mathbb{P}}$-almost surely. The identifications $\ti{g}^{h} = M_{h}(\ut^{h})\ti{p}_{x}^{h}$ and $\ti{g} = \ut^{2} \ti{p}_{x}$ can be achieved as in Lemma 5.6 and Lemma 5.9 in \cite{FischerGruen}. Thus, we have 
	\begin{align}
		M_{h}(\ut^{h})\ti{p}^{h}_{x} \rightharpoonup \ut^{2} \ti{p}_{x}.
	\end{align}
	By the lower semi-continuity of the $L^{2}$-Norm w.r.t. weak convergence and Fatou's lemma, we have
	\begin{align}\notag
		\Ew{\int_{0}^{T} \int_{\cO} \left( \ut^{2}\ti{p}_{x} \right)^{2} \dx \ds}
		&\le 
		\Ew{\underset{h \rightarrow 0}{\liminf}\int_{0}^{T} \underset{x\in \cO}{\sup}((u^{h})^{2}) \int_{\cO} M_{h}(\ut^{h})(\ti{p}^{h}_{x})^{2} \dx \ds	}
		\\ \notag
		&\le
		\underset{h \rightarrow 0}{\liminf}\, 
		\Ew{\int_{0}^{T} \underset{x\in \cO}{\sup}((u^{h})^{2}) \int_{\cO} M_{h}(\ut^{h})(\ti{p}^{h}_{x})^{2} \dx \ds	}
		\\
		&\le C(\e, \un,T) 
		< \infty \, ,
	\end{align}
	where we could use \eqref{lem:UniformEstimate-u^2ph4} due to \eqref{lem:UniformEstmate-u^2ph2}. 
	The arguments to show \eqref{lem:UniformEstimate-u^2ph2} and \eqref{lem:UniformEstimate-u^2ph3} are similar.
	\proofend
	
\section{A Combined $\alpha$-Entropy-Energy Estimate}
	\label{sec:Energy-alpha-Entrop}
	In this section, we prove a new estimate satisfied by solutions to equation \eqref{eq:FischerGruenModified} which is independent of $\e>0.$ This estimate will be the key to pass to the limit $\e\to 0$ and this way to establish the existence of martingale solutions to equation \eqref{eq:StfeStrato}.
	Abbreviating $F(u):=u^{-p}$, equation \eqref{eq:FischerGruenModified} can equivalently be written
	
	\begin{align}\label{eq:stfe-Strat-equ} 
		\dif \ue =  - ((\ue)^2(\uexx - \e \W'(\ue))_{x})_{x} \dt  
		+ 
		\Cstr \uexx \dt
		+ \sum_{k \in \Z} \lal (\gl\ue)_{x} \dble \, .  
	\end{align}	 
	Theorem \ref{theo:existence-stfe-stratonovich} guarantees the existence of a family  $(\ue)_{\e \in (0,1)}$ of solutions to \eqref{eq:stfe-Strat-equ}.
	\begin{the1}\label{theo:ExistenceApproxSolutions}
		Let $\e \in (0,1)$ and  $\ue$ be a solution to  \eqref{eq:stfe-Strat-equ}. Then, for $q \geq 1$ and $\alpha \in (-\frac{1}{3},0)$ the $\e$-independent estimate
		\begin{align}\notag  \label{est:ApproxSolutions}
			\mathbb{E}   &\left[  \sup_{t\in [0,T]} \left (\int_{\cO}\frac{1}{2}\lvert \ue_{x}\rvert^{2} \right.	+ \e\W(\ue) \dx  \right)^{q}			  
			+
			\sup_{t\in [0,T]} \left( \int_{\cO}\frac{1}{\alpha (\alpha+1)} \ueap- \frac{1}{\alpha}\ue + 
			\frac{1}{\alpha +1} \dx \right)^{q}   \\ \notag
			+
			&\left( \intsptiT{(\ue)^{2}(\uexx-\e \W'(\ue))_{x}^{2}}\right)^{q} 
			\\ \notag
			+
			&\left(\intsptiT{((\ue)^{\frac{\alpha + 3}{4}})_{x}^{4}}\right)^{q}
			+
			\left(\intsptiT{((\ue)^{\frac{\alpha+3}{2}})^{2}_{xx}}\right)^{q} 
			\\ 
			+
			&\left.\left(\intsptiT{\ueap (\uex)^{2}  \e\W''(\ue)}\right)^{q} \right]  
			\leq C(T,q,\un).
		\end{align}
		holds.
	\end{the1}
	Our strategy to prove Theorem \ref{theo:ExistenceApproxSolutions} is to combine It\^o's formula and a stopping time argument. 
	Let us for $\e > 0$ consider energies 
	\begin{align}\label{def:Energy}
		\mathcal{E}_{\e}(u) := \frac{1}{2}\int_{\cO}  \betr{\ux}^{2} \dx +  \int_{\cO} \Pi_\e(u)	 \dx
	\end{align}
	as well as  stopping times 
	\begin{align}\label{def:stoppingTimes}
		T_{\sigma}:= T \wedge \inf\{t \in [0,T] \vert \mathcal{E}_{\e}(\ue) \ge \sigma^{-1} 	\}
	\end{align}
	for positive parameters $\sigma$.
	We further introduce the following cut-off versions of solutions $\ue$, where we skip the index $\e$:
	\begin{align}\label{def:ModifiedSolutions}
		\us(\cdot, t) :=
		\begin{cases}
			\ue(\cdot, t) \quad & t \in [0, T_{\sigma}]
			\\
			\ue(\cdot, T_{\sigma}) \quad & t \in (T_{\sigma},\Tmax] \, .
		\end{cases}
	\end{align}	
	Moreover, we set $\ps := -\usxx + \e \W'(\us)$.	
	\begin{lem1}\label{lem:alpha-Entropy-Energy}
		For $\alpha \in (-\frac{1}{3},0)$, $q \ge 1$, and a constant $C(T,q, \un)$ that is independent of $\e$ we have the estimate
		\begin{align}\label{eq:alpha-Entropy-Energy} \notag
			&\Ew{ \sup_{t\in [0,\Tmax]} 
				\left( 
				\frac{1}{2}\int_{\cO} (\us)_{x}^{2}
				+
				\e \W(\us)
				\dx
				\right)^{q} 
				+
				\sup_{t\in [0,\Tmax]}
				\left(
				\int_{\cO} 
				\frac{1}{\alpha(\alpha +1)} \us^{\alpha+1} - \frac{1}{\alpha}\us + \frac{1}{\alpha+1} \dx 
				\right)^{q}
			}
			\\ \notag
			+
			&\Ew{ \left(
				\int_{0}^{t\wedge T_{\sigma}} \int_{\cO}
				(\us)^{2} \psx^{2} \dx \ds
				\right)^{q}			
			}
			+
			\Ew{\left(
				\int_{0}^{t\wedge T_{\sigma}} \int_{\cO}
				(\us)_{xx}^{2} (\us)^{\alpha +1} \dx \ds
				\right)^{q}
			}
			\\ \notag
			+
			&\Ew{\left(
				\int_{0}^{t\wedge T_{\sigma}} \int_{\cO}
				\frac{\abs{\alpha(\alpha+1)}}{3} \us^{\alpha-1} (\us)_{x}^{4} \dx \ds
				\right)^{q}
			}
			\\ 
			+
			&\Ew{\left(
				\int_{0}^{t\wedge T_{\sigma}} \int_{\cO}
				\us^{\alpha+1} (\us)^{2}_{x} \e\W''(\us)  \dx \ds
				\right)^{q}
			} 
			\le
			C(T,q, \un) \, .
		\end{align}
	\end{lem1}
	The proof of Lemma \ref{lem:alpha-Entropy-Energy} will be given below. Let us first assume it to hold and see how \eqref{est:ApproxSolutions} can be derived from it. We need
	\begin{lem1}\label{lem:ConvergenceStoppingTimes}
		We have $\Liml_{\sigma \rightarrow 0} T_{\sigma} = \Tmax$ \, $\mathbb{P}^{\e}$-almost surely.
	\end{lem1}
	\proof
	By estimate \eqref{est:EnEntropFischeGruen} we have for arbitrary, but fixed, $\e \in (0,1)$ and $q \ge 1$
	\begin{align}\notag
		\Ew{
			\supl_{t \in [0,\Tmax]} \mathcal{E}_{\e}(\ue)^{q} 
		}
		\le C(\e) <
		\infty \, .
	\end{align}
	This implies
	\begin{align}\notag
		\mathbb{P}^{\e}
		\left(
		\{T_{\sigma} < \tau\} 
		\right)
		=
		\mathbb{P}^{\e}
		\left(
		\left\lbrace \omega \in \Omega \biggr\vert \supl_{t \in [0,\tau]} \mathcal{E}(\ue(\cdot,t)) \ge \sigma^{-1} \right\rbrace   
		\right) 
		\le
		C(\e) \sigma
	\end{align}
	for any $\tau \in (0,\Tmax]$.
	Hence, for $\sigma \rightarrow 0$ we have $T_{\sigma}\rightarrow \Tmax$ in probability and therefore $T_{\sigma } \rightarrow \Tmax$ $\mathbb{P}^{\e}$-almost surely for a subsequence.
	\proofend
	
	\begin{myproof}{ of Theorem \ref{theo:ExistenceApproxSolutions}}
		From Lemma~\ref{lem:Minimum-ue} we infer that $\ue$ is strictly positive on $\ort\times[0,T]$ $\probab^{\e}$-almost surely. In combination with Lemma~\ref{lem:ConvergenceStoppingTimes}, we find that the sets $A_\sigma:=\{\omega\in\Omega:\ue(\cdot,\omega)\geq \sigma \text{ on } \ort\times[0,T]\}$ tend for $\sigma\to 0$ to $\Omega$ up to a set of measure zero. Hence, using nonnegativity of the terms on the left hand side of \eqref{eq:alpha-Entropy-Energy} and monotone convergence, Theorem~\ref{theo:ExistenceApproxSolutions} is proven.
	\end{myproof}
	
	The rest of this section is devoted to the proof of Lemma \ref{lem:alpha-Entropy-Energy}. 
	\begin{myproof}{ of Lemma \ref{lem:alpha-Entropy-Energy}}
		Consider the operators
		\begin{align}
			&E_{1}: u \mapsto \frac{1}{2}\int_{\cO} u^{2} \dx \, ,
			\\ 
			&E_{2}: u \mapsto \e \int_{\cO} \W (\eta(u)) \dx \, ,
			\\ 
			&\Gmcas: u \mapsto \int_{\cO} G_{\alpha}(\eta(u)) \dx \, ,
		\end{align}
		where $\eta$ is a positive smooth cut-off function corresponding to the lower bound of $\ue$ provided by Lemma \ref{lem:Minimum-ue}, cf. also \eqref{def:CutoffFunc}. Moreover, $G_{\alpha}$ is a standard $\alpha$-entropy used for the thin-film equation, i.e.
		\begin{align}
			G_{\alpha}(u) = \frac{1}{\alpha(\alpha +1)} u^{\alpha+1} - \frac{1}{\alpha}u + \frac{1}{\alpha+1} 
			> 0 \, .
		\end{align} 
		It\^o's formula, applied separately for each of these operators, see appendix \ref{appendix:Itoformula} for details, yields 
		
		\begin{align}\label{eq:CombinedEnEntr}
			\notag 
			&\int_{\cO} \frac{1}{2}(\us(t))_ {x}^{2} + \e \W(\us(t))\dx
			+
			\Gmcas (\us(t)) 
			\\ \notag
			&+
			\int_{0}^{t\wedge T_{\sigma}} \int_{\cO}
			(\us)^{2} \psx^{2}   
			+
			\int_{0}^{t\wedge T_{\sigma}} \int_{\cO}
			C_{Strat} (\us)_{xx}^{2} 
			\dx\ds
			\dx \ds
			\\ \notag
			&+
			\int_{0}^{t\wedge T_{\sigma}} \int_{\cO}
			C_{Strat} \usx^{2} \e \W''(\us)
			\dx \ds
			\\ \notag
			=
			&E_{1}((\uen)_{x}) +E_{2}(\uen) + \Gmcas (\uen) 
			\\ \notag
			&+
			\sum_{ k \in \Z} 
			\int_{0}^{t\wedge T_{\sigma}} \int_{\cO} 
			(\us)_{x}\lal (\us \gl)_{xx}
			+
			\lal \e \W'(\us) (\us \gl)_{x} 
			\dx \dbl
			\\ \notag
			&+
			\frac{1}{2} \int_{0}^{t\wedge T_{\sigma}} \int_{\cO} \sum_{ k \in \Z} \lal^{2} (\us \gl)^{2}_{xx} \dx\ds
			\\ \notag
			&+
			\frac{1}{2} \int_{0}^{t\wedge T_{\sigma}} \int_{\cO}
			\sum_{ k \in \Z} \lal^{2} \e \W''(\us)(\us \gl)^{2}_{x} \dx \ds 
			\\ \notag
			&+
			\sum_{ k \in \Z} 
			\int_{0}^{t\wedge T_{\sigma}} \int_{\cO} 
			\lal
			\left(\frac{1}{\alpha} \us^{\alpha} - \frac{1}{\alpha}\right) (\us \gl )_{x}
			\dx \dbl 
			\\ \notag
			&+
			\int_{0}^{t\wedge T_{\sigma}} \int_{\cO}
			(-(\us)^{2} \psx - C_{Strat} (\us)_{x})
			\left(\frac{1}{\alpha} \us^{\alpha} - \frac{1}{\alpha}
			\right)_{x} 
			\dx \ds
			\\ \notag
			&+	
			\frac{1}{2} \int_{0}^{t\wedge T_{\sigma}} \int_{\cO} \sum_{ k \in \Z} \lal^{2} \us^{\alpha-1} (\us \gl)^{2}_{x} \dx \ds 
			\\
			:=
			&E_{1}((\uen)_{x}) +E_{2}(\uen) + \Gmcas (\uen) 
			+
			R_{1} + \dots + R_{6} \, .
		\end{align}

		Let us now derive an estimate of equation \eqref{eq:CombinedEnEntr}. 
		We will frequently use the relations \eqref{eq:ONBRelations1}-\eqref{eq:ONBRelations6} and \eqref{Lambda-Wachstumsbedingung}. Choosing $\sigma$ small enough, we may assume that $\eta(\uen) = \uen$. 
		Then, since $\Lambda^{\e} = \mathbb{P} \circ (\un + \e^{\theta})^{-1}$, (H2), $\theta \in (0,\frac{1}{p})$, and $\W(x) = x^{-p}$, we get
		\begin{align}\notag
			\mathbb{E}^{\e} \left[ E_{1}((\uen)_{x}) + E_{2}(\uen) \right]
			&=
			\mathbb{E}^{0} \left[\frac{1}{2} \int_{\cO} (\un +\e^{\theta})^{2}_{x}  
				+
				\e (\un + \e^{\theta})^{-p} \dx \right]
			\\ \notag
			&\le
			\mathbb{E}^{0} \left[ C \int_{\cO} (\un)_{x}^{2} + \e \e^{-\theta p} \dx \right]
			\le C
		\end{align}
		and for $\alpha\in (-1,0)$ 
		\begin{align}
			\mathbb{E}^{\e} \left[ \Gmcas(\uen) \right]
			&=
			\mathbb{E}^{0} \left[ \int_{\cO} \frac{1}{\alpha(\alpha+1)} (\un + \e^{\theta})^{\alpha+1} 
				-
				\frac{1}{\alpha} (\un +\e^{\theta})
				+
				\frac{1}{\alpha +1} \dx \right]
			\le
			C \, .
		\end{align}
		Here, the superscript indicates that the expectation is computed with respect to $\mathbb{P}^{\e}$ and $\mathbb{P}$, respectively. \newline
		Ad $R_{2}$:
		We have
		\begin{align} \notag
			&\frac{1}{2} \int_{0}^{t\wedge T_{\sigma}} \int_{\cO} \sum_{ k \in \Z} \lal^{2} (\us \gl)^{2}_{xx} \dx \ds  
			\\ \notag
			=
			&\frac{1}{2} \int_{0}^{t\wedge T_{\sigma}} \int_{\cO}
			\sum_{k = 1}^{\infty} k^{4} \lal^{2} \frac{32 \pi^{4}}{L^{5}} 		
			\us^{2}
			+
			2 \sum_{k = 1}^{\infty}  k^{2} \lal^{2}\frac{8\pi^{2}}{L^{3}}
			(\us)^{2}_{x}
			+
			4 \sum_{k = 1}^{\infty} k^{2} \lal^{2}\frac{8\pi^{2}}{L^{3}} 
			(\us)^{2}_{x}
			\\ \notag
			&\qquad \qquad \quad+
			\left(
			\frac{\lan^{2}}{L} 
			+
			\sum_{k = 1}^{\infty} \frac{2\lal^{2}}{L}	
			\right)
			(\us)_{xx}^{2}
			\dx \ds
			\\ 
			:=
			&A +B +C +D \, .
		\end{align}
		By means of Poincar\'e's inequality, A, B, and C will become Gronwall terms, while D also appears on the left hand side of equation \eqref{eq:CombinedEnEntr} and therefore cancels out.
		
		Ad $R_{3}$:
		\begin{align}\notag
			&\frac{1}{2}  \int_{0}^{t\wedge T_{\sigma}} \int_{\cO}
			\sum_{k \in \Z}
			\lal^{2} \e \W''(\us) 
			\left(
			(\gl)^{2}_{x} (\us)^{2} + 2 \glx(\gl) \us (\us)_{x}
			+ 
			(\gl)^{2} (\us)^{2}_{x}
			\right)
			\dx \ds
			\\ 
			&\frac{1}{2}
			\int_{0}^{t\wedge T_{\sigma}} \int_{\cO} 
			\sum_{k = 1}^{\infty}
			k^{2} \lal^{2}\frac{8\pi^{2}}{L^{3}}
			\e p (p+1) \W(\us) \dx\ds
			+
			\int_{0}^{t\wedge T_{\sigma}} \int_{\cO} 
			\e \, 
			C_{Strat}
			\W''(\us) (\us)_{x}^{2} 
			\dx \ds\, ,
		\end{align}
		where the first term will be a Gronwall term, the second cancels out against the corresponding term on the left hand side of \eqref{eq:CombinedEnEntr}.
		
		Let us now discuss the contributions of the entropy.
		\newline
		Ad $R_{5}$: 
		\begin{align} \notag
			&\int_{0}^{t\wedge T_{\sigma}} \int_{\cO}
			(-(\us)^{2} \psx - C_{Strat} (\us)_{x})
			\left(
			\frac{1}{\alpha} \us^{\alpha} - \frac{1}{\alpha}
			\right)_{x} \dx \ds
			\\ \notag
			=
			&\int_{0}^{t\wedge T_{\sigma}} \int_{\cO}
			-(\us)^{2}_{xx} \us^{\alpha+1} 
			+ \frac{\alpha (\alpha+1)}{3} \us^{\alpha-1} (\us)_{x}^{4} 
			-\e \W''(\us) \us^{\alpha+1} (\us)_{x}^{2}
			\dx \ds
			\\ 
			&-
			\int_{0}^{t\wedge T_{\sigma}} \int_{\cO}
			C_{\text{Strat}} \us^{\alpha -1} (\us)^{2}_{x}
			\dx \ds \, .
		\end{align}
		Since $\alpha (\alpha+1) < 0$, the first integral is a good term while the second one will cancel out as the following calculation shows.
		
		Ad $R_{6}$:
		\begin{align}\notag
			&\frac{1}{2} 
			\int_{0}^{t\wedge T_{\sigma}} \int_{\cO} 
			\sum_{ k \in \Z} \lal^{2} \us^{\alpha-1} (\us \gl)^{2}_{x} \dx \ds 
			\\ 
			=
			&\frac{1}{2} \int_{0}^{t\wedge T_{\sigma}} \int_{\cO} \sum_{ k = 1}^{\infty}  
			\lal^{2} k^{2}\frac{8\pi^{2}}{L^{3}}
			\us^{\alpha+1} 
			\dx\ds
			+	
			\int_{0}^{t\wedge T_{\sigma}} \int_{\cO}
			\frac{1}{2}
			\left(
			\frac{\lan^{2}}{L} 
			+
			\sum_{k = 1}^{\infty} \frac{2\lal^{2}}{L}
			\right)
			\us^{\alpha-1} (\us)_{x}^{2} \dx \ds \, .
		\end{align}
		Here the first term can be estimated by means of Young's and Poincar\'e's inequalities and serve as a Gronwall term.
		
		Collecting all terms, rearranging, and combining the constants,
		as well as applying $q$'th powers, $q \ge 1$, suprema, and expectation, we get with $E(v) := E_{1}(v_{x}) + E_{2}(v)$ 
		\begin{align}\notag
			&\Ew{ \sup_{t\in [0,\Tmax]} 
				E(\us(t))^{q} 
				+
				\sup_{t\in [0,\Tmax]}
				\Gmcas(\us(t))^{q}
			}
			\\ \notag
			+
			&\Ew{ \left(
				\int_{0}^{t\wedge T_{\sigma}} \int_{\cO}
				(\us)^{2} \psx^{2} \dx \ds
				\right)^{q}			
			}
			+
			\Ew{\left(
				\int_{0}^{t\wedge T_{\sigma}} \int_{\cO}
				(\us)_{xx}^{2} (\us)^{\alpha +1} \dx \ds
				\right)^{q}
			}
			\\ \notag
			+
			&\Ew{\left(
				\int_{0}^{t\wedge T_{\sigma}} \int_{\cO}
				\frac{\abs{\alpha(\alpha+1)}}{3} \us^{\alpha-1} (\us)_{x}^{4} \dx \ds
				\right)^{q}
			}
			\\ \notag
			+
			&\Ew{\left(
				 \int_{0}^{t\wedge T_{\sigma}} \int_{\cO}
				\e \W''(\us) \us^{\alpha+1} (\us)^{2}_{x} \dx \ds
				\right)^{q}
			}
			\\ 
			\le
			&C^{q} 
			+
			C
			\Ew{
				\int_{0}^{t\wedge T_{\sigma}} \int_{\cO} E(\us)^{q} 
				\dx \ds
			}
			+ \Ew{	
				\sup_{t\in [0,\Tmax]}\abs{R_{1}} ^{q}
			}
			+
			\Ew{	
				\sup_{t\in [0,\Tmax]} \abs{R_{4}}^{q}
			} \, .
		\end{align}		
		In order to use the Burkholder-Davis-Gundy inequality, we first consider for $s \in [0,T]$ the 
		operator $\tau_{1}(s): Q^{\frac{1}{2}} L^{2}(\cO) \rightarrow \R$ with
		\begin{align}\label{BDG-Process1} \notag
			\tau_{1}(s)(v) &:= \chi_{T_{\s}}(s)\int_{\cO} (\us(s))_{x} ( \us(s) \sum_{ i \in \Z} \skp{g_{i}}{v}g_{i} )_{xx}
			\\
			& \quad +
			\e \W'(\us(s))(\us(s)\sum_{ i \in \Z} \skp{g_{i}}{v}g_{i})_{x} \dx \, .
		\end{align}
		Here we set $\chi_{T_{\s}} := \chi_{[0,T_{\s}]}$.
		Let us estimate the Hilbert-Schmidt norm of \eqref{BDG-Process1}. For the proof, $\delta$ and $\ti{\delta}$ are positive constants arising from Young's inequality. For better readability, we will skip the argument $s$. 
		\begin{align}\notag
			&\Norm{\tau_{1}(s)}^{2}_{L_{2}(Q^{1/2}\lpr{2}{\cO},\R)} 
			=
			\sum_{ k \in \Z} \abs{\tau_{1}(s) (\gl \lal)}^{2}
			\\ \notag
			=
			&\sum_{ k \in \Z} \chi_{T_{\s}}
			\abs{
				\int_{\cO}
				(\us)_{xx} \usx \lal \gl 
				+ \usxx \us \lal \glx 
				+ \e p (p+1) \us^{-p-1} \usx \lal \gl 
				\dx 
			}^{2}
			\\ \notag
			=
			& \chi_{T_{\s}} C
			\left(
			\sum_{ k \in \Z} \lal^{2} k^{2} 
			\left(
			\int_{\cO}
			\frac{1}{2} \usx^{2}\dx 
			\right)^{2}
			+
			\sum_{ k \in \Z} \lal^{2} k^{2} 
			\int_{\cO}
			\usxx^{2} \us^{2} 
			\dx
			\right.
			\\ \notag
			&\qquad+
			\left.
			\sum_{ k \in \Z} \lal^{2} k^{2}  (p+1)^{2}
			\left(
			\e \int_{\cO}
			\W(\us) 
			\dx
			\right)^{2}
			\right)			
			\\ 
			\le
			&\chi_{T_{\s}} 
			\left(
			C(\lal, p) E(\us)^{2} 
			+
			C(\lal) \Norm{\us^{-\alpha+1}}_{\infty} 
			\int_{\cO} \usxx^{2} \us^{\alpha+1} \dx
			\right)	\, .
		\end{align}
		Hence, 
		\begin{align}\notag
			&\qquad \Ew{\sup_{t\in [0,\Tmax]}
				\abs{R_{1}}^{q}
			}
			\le
			C 
			\Ew{
				\left(
				\int_{0}^{t \wedge T_{\s}}
				%C 
				E(\us)^{2} + 
				%C
				\Norm{\us^{-\alpha+1}}_{\infty} 
				\int_{\cO} 
				\usxx^{2} \us^{\alpha+1} 
				\dx \ds
				\right)^{\frac{q}{2}}
			}
			\\ \notag
			&\le
			C(q)
			\Ew{ C
				\int_{0}^{T_{\sigma}}
				E(\us)^{q} 
				\ds
			}
			+
			C(p,\delta,T) \, \ti\delta \, 
			\Ew{
				\sup_{t\in [0,T_{\sigma}]}
				\Norm{\us}_{\Hsob{1}{\cO}}^{2q} 
			}
			+
			\ti{\delta}^{-1} C(q,\delta,T,L) 
			\\
			&\quad+ C(q)  \delta
			\Ew{
				\left(
				\int_{0}^{t \wedge T_{\s}}
				\int_{\cO} 
				\usxx^{2} \us^{\alpha+1} 
				\dx \ds
				\right)^{q}
			} \, .
		\end{align}
		The first term in the last estimate is a Gronwall term while the others can be absorbed.
		We proceed to do the same for the stochastic integral $R_{4}$. 
		The corresponding operator $\tau_{2}(s): Q^{\frac{1}{2}}L^{2}(\cO) \rightarrow \R$ now reads 
		\begin{align}\label{BDG-Process2}\notag
			\tau_{2}(s)(v) := \chi_{T_{\s}}(s) 
			&\int_{\cO} 
			\left(
			\frac{1}{\alpha} \us^{\alpha}(s) - \frac{1}{\alpha}
			\right)
			(\us(s) \sum_{ i \in \Z} \skp{g_{i}}{v}g_{i})_{x} 
			\dx
		\end{align}
		for $s \in [0, \Tmax]$. It follows
		\begin{align}\notag
			\Norm{\tau_{2}(s)}^{2}_{L_{2}(Q^{1/2}\lpr{2}{\cO},\R)} 
			=
			&\sum_{ k \in \Z} \chi_{T_{\s}}
			\abs{ \int_{\cO} 
				\left(
				\frac{1}{\alpha} \us^{\alpha} - \frac{1}{\alpha}
				\right)
				(\us \sum_{ i \in \Z} \skp{\gl \lal}{g_{i} } g_{i})_{x}
				\dx
			}^{2}
			\\ \notag
			\le
			&  \sum_{ k \in \Z} \chi_{T_{\s}}
			\int_{\cO}
			\us^{2\alpha} \usx^{2} \gl^{2} \lal^{2}
			\dx
			\\ \notag 
			\le
			&  \sum_{ k \in \Z} \chi_{T_{\s}}
			\lal^{2} 
			\Norm{\gl}_{\infty}^{2}
			\Norm{\us^{\frac{3\alpha+1}{2}}}_{\infty}^{2}
			\int_{\cO}
			\us^{\frac{\alpha-1}{2}} \usx^{2}
			\dx
			\\ 
			\le
			& C(\lal) \chi_{T_{\s}} 
			\Norm{\us^{3\alpha+1}}_{\infty}
			\int_{\cO}
			\us^{\frac{\alpha-1}{2}} \usx^{2}
			\dx\, ,
		\end{align}
		where we used $\alpha \in (-\frac{1}{3},0)$ 
		, the boundedness of $\gl$, and \eqref{Lambda-Wachstumsbedingung}. Then we get
		\begin{align}\notag
			\Ew{\sup_{t\in [0,\Tmax]}
				\abs{R_{4}}^{q}
			}
			&\le
			C	
			\Ew{\left(
				\int_{0}^{t \wedge T_{\s}} 
				C(\lal) 
				\Norm{\us^{3\alpha+1}}_{\infty}
				\int_{\cO}
				\us^{\frac{\alpha-1}{2}} \usx^{2}
				\dx \ds
				\right)^{\frac{q}{2}}		
			}
			\\\notag 
			&\le
			C(q,\lal,T) \, \ti{\delta} \,
			\Ew{  
				\sup_{t\in [0,T_{\sigma}]}  \Norm{\us}_{\Hsob{1}{\cO}}^{2q}
			}
			+
			\ti{\delta}^{-1}
			C(q,\lal,T,L)^{q}
			\\ 
			& \quad+
			C(q) \delta
			\Ew{
				\left(
				\int_{0}^{t \wedge T_{\s}} \int_{\cO}
				\us^ {\alpha-1} \usx^{4} \dx \ds
				\right)^{q} 
			}
			+ \delta^{-q}C(q,T,L)	\, .
		\end{align}
		For the last estimate we have once more used Poincar\'e's inequality.
		The first and the third term can be absorbed while the others are independent of $\e$.
		After absorption, we therefore obtain
		\begin{align}\notag
			&\Ew{ \sup_{t\in [0,\Tmax]} 
				E(\us(t))^{q} 
			}
			+
			\Ew{
				\sup_{t\in [0,\Tmax]}
				\Gmcas(\us(t))^{q}
			}
			\\ \notag
			+
			&\Ew{ \left(
				\int_{0}^{t\wedge T_{\sigma}} \int_{\cO}
				(\us)^{2} \psx^{2} \dx \ds
				\right)^{q}			
			}
			+
			\Ew{\left(
				\int_{0}^{t\wedge T_{\sigma}} \int_{\cO}
				(\us)_{xx}^{2} (\us)^{\alpha +1} \dx \ds
				\right)^{q}
			}
			\\ \notag
			+
			&\Ew{\left(
				\int_{0}^{t\wedge T_{\sigma}} \int_{\cO}
				\frac{\abs{\alpha(\alpha+1)}}{3} \us^{\alpha-1} (\us)_{x}^{4} \dx \ds
				\right)^{q}
			}
			\\ \notag
			+
			&\Ew{\left(
				 \int_{0}^{t\wedge T_{\sigma}} \int_{\cO}
				\e \W''(\us) \us^{\alpha+1} (\us)^{2}_{x} \dx \ds
				\right)^{q}
			}
			\\ 
			\le
			&C(q,\lal,T) 
			+
			C
			\Ew{
				\int_{0}^{T_{\sigma}} E(\us)^{q} \dx \ds	
			} \, .
		\end{align}
		The very last term can be estimated, since for  $t \in [0,\Tmax]$ 
		\begin{align}\label{est:gronwall}
			\Ew{
				E(\us(t))^{q} 	
			} 
			\le
			\Ew{ \sup_{t\in [0,T_{\sigma}]} 
				E(\us(t))^{q} 
			}
			\le 
			C +C \int_{0}^{T_{\sigma}} 
			\Ew{
				E(\us)^{q} 	
			}\ds.
		\end{align}
		We have used Fubini's theorem, which is possible due to the nonnegativity of the energy.
		Gronwall's Lemma shows that the last term in \eqref{est:gronwall} is finite, so we have finished the proof.		
	\end{myproof}	
	
\section{Convergence of Approximate Solutions}
	\label{sec:Jakubowski}
In this section, we pass to the limit $\e\to 0$ with the approximate solutions $u^\e $ and finally, we prove  the main results of the paper, i.e. Theorem~\ref{theo:MainResult} and Corollary~\ref{cor:contactangle}.	
\subsection{Application of Jakubowski's Theorem}
	\label{subsec:Kompaktheit}
	In this subsection we will apply \newline Jakubowski's theorem, cf. \cite{Jakubowski1994}.
	
	We define
	$
	\ve := ((\ue)^{\frac{\alpha + 3}{4}})_{x}
	$,
	$
	\ze := ((\ue)^{\frac{\alpha+3}{2}})_{xx} 
	$, 
	and for $\ti{\gamma} < \gamma$ and $\ti{\sigma} < \sigma$
	\begin{align*}
		\cX_{u} &:=  C^{\tilde{\gamma}, \tilde{\sigma}}(\Otb) 
		\\
		\cX_{\ux} &:=  L^{2}(0,T;L^{2}(\cO))_{weak} 
		\\
		\cX_{v} &:=  L^{4}(0,T;L^{4}(\cO))_{weak} 
		\\
		\cX_{z} &:=  L^{2}(0,T;L^{2}(\cO))_{weak} 
		\\
		\cX_{W} &:=  C([0,T];L^{2}(\cO)) 
		\\
		\cX_{\un} &:= H^{1}_{\text{per}}(\cO).
	\end{align*} 
	Moreover, $\cX := \cX_{u} \times \cX_{\ux} \times \cX_{v} \times \cX_{z} \times \cX_{W} \times \cX_{\un}$.
	\begin{sat1}\label{sat:Strafheit}
		The laws  $\mu_{\ue},\mu_{\uex},\mu_{\ve},\mu_{\ze},\mu_{W^{\e}}, \mu_{\uen}$ of the corresponding random variables are tight on the path spaces $\cX_{u},\cX_{\ux}, \cX_{v}, \cX_{z}, \X_{W}$, and $\cX_{\un}$, respectively.
	\end{sat1}
	
	\proof
	Since $\bar{B}_{R} \subset C^{\gamma ,\sigma}(\Otb)$ is a compact subset of $C^{\ti{\gamma} ,\ti{\sigma}}(\Otb)$, we can use the uniform bound of $\ue $ in $L^{2}(\Omega^{\e};C^{\gamma ,\sigma}(\Otb))$ that follows from Theorem \ref{theo:ExistenceApproxSolutions}. Combining this with Markov's inequality, we get for arbitrary $R > 0$ 
	\begin{align} \notag
		\mu_{\ue}&(C^{\ti{\gamma} ,\ti{\sigma}}(\Otb) \backslash \bar{B}_{R}) \notag
		= \mathbb{P}^{\e} \left[ \Norm{\ue}_{C^{\gamma ,\sigma}(\Otb)} > R \right] 
		\leq
		\frac{\Ew{\Norm{\ue}_{C^{\gamma ,\sigma}(\Otb)}^{2} } }{R^{2}}
		\le  \frac{C}{R^{2}}.
	\end{align}
	Since $C$ is independent of $\e$, for $R$ big enough we get tightness by choosing the compact set $K = \bar{B}_{R}$.
	To prove tightness of $(\mu_{\uex})_{\e \in (0,1)}$, $(\mu_{\ve})_{\e \in (0,1)}$, and $(\mu_{\ze})_{\e \in (0,1)}$, we can use similar arguments and the corresponding terms in the $\alpha$-entropy-energy estimate \eqref{est:ApproxSolutions}.
	The processes $W^{\e}$ share for all $\e \in (0,1)$ their law as they are Q-Wiener processes.
	On the Polish space $C([0,T];L^{2}(\cO))$ the laws  $\mu_{W^{\e}}$ are regular. In particular $\mu_{W^{\e}}$ is regular from the interior, which means
	\begin{align}\notag
		\mu_{W^{\e}}(A) = \sup \{\mu_{W^{\e}}(K): \, K \subset A,\, K \,\, \text{kompakt}	\}
	\end{align}
	for all sets $A$ in the Borel $\sigma$-algebra on $C([0,T];L^{2}(\cO))$. This shows tightness of \newline $(\mu_{W^{\e}})_{\e \in (0,1)}$ on $C([0,T];L^{2}(\cO))$. 
	By the same reasoning, we get tightness of the law of $\un$ on $H^{1}_{per}(\cO) $, i.e. for arbitrary $\delta > 0$, we find a compact set $K_{\delta} \subset H^{1}_{per}(\cO)$ with
	\begin{align}\notag
		\mathbb{P} \left[ \un \in (K_{\delta})^{c} \right] \le \delta .
	\end{align}
	Consider the set $\ti{K_{\delta}} = \{ f + \eta \colon f \in K_{\delta}, \, \eta  \in [0,1]	\}$. Then $\ti{K_{\delta}}$ is compact and $\un \notin K_{\delta} \Leftrightarrow \un + \e^{\theta} \notin \ti{K}_{\delta}$. Hence,
	\begin{align}\notag
		\mathbb{P}^{\e} \left[ \uen \in (\ti{K_{\delta}})^{c} \right] 
		= 	\mathbb{P} \left[ \un + \e^{\theta} \in (\ti{K_{\delta}})^{c} \right] 
		= 	\mathbb{P} \left[ \un \in ({K_{\delta}})^{c} \right] \le \delta \, .
	\end{align}
	This completes the proof. 
	\proofend

	\begin{sat1}\label{sat:JakubowskiAnwendung}
		For subsequences of $\ue, \uex, \ve , \ze, \uen$, and $W^{\e}$, there exist a probability space $(\ti{\Omega}, \ti{\mathcal{F}}, \ti{\mathbb{P}})$, a sequence of $L^{2}(\cO)$-valued stochastic
		processes $\ti{W}^{\e}$ on $(\ti{\Omega}, \ti{\mathcal{F}}, \ti{\mathbb{P}})$, sequences of random variables
		\begin{align*}
			\uet: &\ti{\Omega} \rightarrow C^{\ti{\gamma},\ti{\sigma}}(\Otb ), 
			\\
			\wet: &\ti{\Omega} \rightarrow L^{2}(0,T;L^{2}(\cO)), 
			\\
			\vet: &\ti{\Omega} \rightarrow \lpr{4}{0,T;\lpr{4}{\cO}}, 
			\\
			\zet: &\ti{\Omega} \rightarrow \lpr{2}{0,T;\lpr{2}{\cO}}, 
			\\
			\uent: &\ti{\Omega} \rightarrow H^{1}_{\text{per}}(\cO), 
		\end{align*} 	
		random variables
		\begin{align*}
			&\ut \in L^{2}(\ti{\Omega};C^{\ti{\gamma},\ti{\sigma}}(\Otb )), \\
			&\wt \in  L^{2}(\ti{\Omega};L^{2}(0,T;L^{2}(\cO))), 
			\\
			&\vt \in L^{4}(\ti{\Omega};\lpr{4}{0,T;\lpr{4}{\cO}}),	
			\\ 
			&\zt \in L^{2}(\ti{\Omega};\lpr{2}{0,T;\lpr{2}{\cO}}), 
			\\
			&\unt \in L^{2}(\ti{\Omega};H^{1}_{\text{per}}(\cO)) ,
			\\ 
		\end{align*}
		as well as a $L^{2}(\cO)$-valued process $\ti{W}$ such that
		\begin{enumerate}
			\item
			for all $\e \in  (0,1)$  the law of $(\uet,\wet, \vet, \zet,\ti{W}^{\e}, \uent)$ on $\cX$ w.r.t. the measure $\ti{\mathbb{P}}$ equals the law of $(\ue, \uex, \ve, \ze,W^{\e}, \uen)$ w.r.t. $\mathbb{P}^{\e}$.
			\item 
			as $\e \rightarrow 0$, the sequence $(\uet, \wet, \vet, \zet,\ti{W}^{\e}, \uent)$ converges  $\ti{\mathbb{P}}$-almost surely to \newline
			$(\ut, \wt, \vt, \zt,\ti{W}, \unt)$ in the topology of $\cX$.		
		\end{enumerate}
	\end{sat1}
%	\proofend
	Next, we identify the new sequences on $(\ti{\Omega},\ti{\mathcal{F}}, \ti{\mathbb{P}})$. 
	
	\begin{lem1}\label{lem:IdentifikationFormel}
		We have 
		$
		\wet = \uext 
		$
		,
		$
		\vet = ((\uet)^{\frac{\alpha +3}{4}})_{x}
		$,
		as well as
		$
		\zet = ((\uet)^{\frac{\alpha +3}{2}})_{xx} 
		$
		$\mathbb{\ti{P}}$-almost surely.
	\end{lem1}
	\proof
	The mapping
	\begin{align*}
		\lpr{2}{0,T;\lpr{2}{\cO}} \rightarrow  \mathbb{R}, \quad u \mapsto \intsptiT{u}	
	\end{align*}
	is Borel-measurable.
	Hence, for arbitrary $\phi \in C^{\infty}_{c}(\Ot)$ we have due to the equality of laws stated in Theorem \ref{sat:JakubowskiAnwendung}
	\begin{align*}
		&\Ew{\abs{\intsptiT{\wet \phi}+\intsptiT{\uet\phi_{x}}}} 
		\\ \notag
		=&\Ew{\abs{\intsptiT{\uex \phi}-\intsptiT{\ue_{x}\phi }} } 
		\\ 
		=&0 \,	.	 
	\end{align*}
	The other statements follow by similar reasoning.
	\proofend
	
	As in \cite{FischerGruen}, we consider the filtrations $(\ti{F}_{t})_{t \ge 0 }$ and $(\ti{F}^{\e}_{t})_{t \ge 0}$ with
	\begin{align}\label{def:Ft}
		\ti{F}_{t} := \sigma(\sigma(r_{t}\ut,r_{t}\Wt) \cup \{N \in \ti{F}:\ti{\mathbb{P}}(N)=0
		\} \cup \sigma(\unt))
	\end{align}  
	and
	\begin{align}\label{def:F}
		\ti{F}^{\e}_{t} := \sigma(\sigma(r_{t}\uet,r_{t}\Wet) \cup \{N \in \ti{F}:\ti{\mathbb{P}}(N)=0
		\} \cup \sigma(\uent)) \, .
	\end{align}  
	Here, $r_{t}$ is the restriction of a mapping on $[0,T]$ to the time interval $[0,t]$, $t \in [0,T]$.
	The proof of the next lemma can be found in \cite{FischerGruen} Lemma 5.7.
	
	\begin{lem1}\label{lem:WienerProcesses}
		The processes $\Wet$ and $\Wt$ are Q-Wiener processes which are adapted to the filtrations $(\ti{F}^{\e}_{t})_{t \ge 0}$ and $(\ti{F}_{t})_{t \ge 0 }$, respectively. We have
		\begin{align}\label{lem:WienerProcesses-h1}
			\Wet(t) = \sum_{k \in \Z} \lal \blet(t)\gl
		\end{align}
		and
		\begin{align}\label{lem:WienerProcesses-h2}
			\Wt(t) = \sum_{k \in \Z} \lal \blt(t)\gl
		\end{align}
		with families $(\blet)_{k \in \Z}$ and $(\blt)_{k \in \Z}$ of i.i.d. standard Brownian motions w.r.t. $(\ti{F}^{\e}_{t})_{t\ge 0}$ and $(\ti{F}_{t})_{t\ge 0}$, respectively.
	\end{lem1}
	
	For the limits $\wt$, $\vt$ and $\zt$, we get the following identities.
	
	\begin{lem1}\label{lem:Identifikation-wt-vt-zt}
		We have $\mathbb{\ti{P}}$-almost surely
		$
		\wt =  \uxt 
		$,	
		$
		\vt =  (\ut^\frac{\alpha+3}{4})_{x} 
		$, and
		$
		\zt = (\ut^{\frac{\alpha+3}{2}})_{xx} 
		$
		.	
	\end{lem1}   
	\proof 
	Exemplarily we show the first statement. From our convergence results in Theorem \ref{sat:JakubowskiAnwendung} and by integration by parts, we deduce for all test functions $\phi \in C_{c}^{\infty}(\Ot)$ 
	\begin{align}\notag
		\intsptiT{	\wt \phi } \leftarrow \intsptiT{\wet \phi } 
		&= \intsptiT{	\uext \phi }  
		\\ \notag
		&= - \intsptiT{ \uet \phi_{x} } 
		\\ \notag
		&\rightarrow  
		- \intsptiT{\ut \phi_{x} } 
	\end{align}
	$\mathbb{\ti{P}}$-almost surely. This gives the first equality.
	\proofend
	
\subsection{Convergence Results of the Deterministic Terms}
	In the next lemmas we establish convergence of the deterministic terms, corresponding to the weak formulation \eqref{eq:MainResultWeakForm}.
	\begin{lem1}\label{lem:schwacheKonvergenz-uex^3}
		The sequence $\uet$ admits a subsequence such that for a function $\zeta$ on $\ti{\Omega} \times \cO\times (0,T) $ 
		\begin{align}
			(\uext)^{3} \rightharpoonup \zeta
		\end{align} 
		weakly in $\lpr{\frac{4}{3}}{\cO_{T}}$ $\ti{\mathbb{P}}$-almost surely. 
	\end{lem1}
	\proof 
	From Lemma \ref{sat:JakubowskiAnwendung}, Lemma \ref{lem:IdentifikationFormel}, and Lemma  \ref{lem:Identifikation-wt-vt-zt}  in the previous subsection, we know
	$
	((\uet)^{\frac{\alpha+3}{4}})_{x} \rightharpoonup 	(\ut^{\frac{\alpha+3}{4}})_{x}
	$
	$\ti{\mathbb{P}}$-almost surely in $\lpr{4}{0,T;\lpr{4}{\cO}}$. By the identity
	\begin{align}\label{eq:Darstellung-uxe^3}
		(\uext)^{3} = \left(\frac{4}{\alpha+3}\right)^{3} (\uet)^{\frac{3(-\alpha +1)}{4}} 	((\uet)^{\frac{\alpha+3}{4}})^{3}_{x} \,\,  , 
	\end{align}
	the $\ti{\mathbb{P}}$-almost surely uniform boundedness of $\uet$  in $\lpr{\infty}{\cO_{T}}$, and the positivity of $-\alpha +1$ we conclude 
	\begin{align}\notag
		\intsptiT{((\uext)^{3})^{\frac{4}{3}}} 
		&=  
		C	\intsptiT{  (\uet)^{-\alpha +1}	((\uet)^{\frac{\alpha+3}{4}})^{4}_{x}	} 
		\le C
	\end{align}
	$\ti{\mathbb{P}}$-almost surely. The result then follows by the reflexivity of $\lpr{\frac{4}{3}}{\cO}$.
	\proofend
	
	Let $S$ be a set. $L^{p-}(S)$ denotes the space of functions that are contained in every space  $L^{q}(S)$, where $1\le q < p$.
	
	\begin{lem1}\label{lem:starkeKonvergenz-L4-}
		We have $\ti{\mathbb{P}}$-almost surely
		\begin{align}
			((\uet)^{\frac{\alpha+3}{4}})_{x} \rightarrow  ((\ut)^{\frac{\alpha+3}{4}})_{x}
		\end{align} strongly in $\lpr{4-}{[\ut>0]}$.
	\end{lem1}
	\proof
	For arbitrary $p > 1 $, we have the strong convergence of $\uet$ in $ \lpr{p}{\cO_{T}}$ $\ti{\mathbb{P}}$-almost surely, cf. Theorem \ref{sat:JakubowskiAnwendung}. Thus, 
	$
	(\uet)^{\frac{\alpha+3}{2}} \rightarrow \ut^{\frac{\alpha + 3}{2}} 
	$
	in $\lpr{2}{\cO_{T}}$ $\ti{\mathbb{P}}$-almost surely.
	By Riesz' theorem
	\begin{align}\label{lem:starkeKonvergenz-L4-h1}
		\underset{h \rightarrow 0}{\lim} \Norm{(\uet)^{\frac{\alpha+3}{2}}(\cdot , \cdot+h) -(\uet)^{\frac{\alpha+3}{2}}(\cdot , \cdot)}_{L^{2}((0,T-h);L^{1}(\mathbb{\cO}))} = 0 	
	\end{align}
	follows. Furthermore, \eqref{est:ApproxSolutions} implies
	\begin{align}\label{lem:starkeKonvergenz-L4-h2}
		\underset{\e \in (0,1)}{\sup}	\Norm{((\uet)^{\frac{\alpha+3}{2}})_{xx}}_{\lpr{2}{\Ot}} \le C \, . 
	\end{align} 
	Hölder's inequality as well as the uniform bound of $\ue$ in $C^{\ti{\gamma} ,\ti{\sigma}}(\Otb)$ show
	\begin{align}\label{lem:starkeKonvergenz-L4-h3} \notag
		\intsptiT{( (\uet)^{\frac{\alpha + 3}{2}} )^{2}_{x}} 
		= 
		C\intsptiT{(\uext)^{2} (\uet)^{\frac{\alpha -1}{2}} (\uet)^{\frac{\alpha +3}{2}}} 
		\\ 
		\le
		C \left( \intsptiT{(\uext)^{4} (\uet)^{\alpha -1} } \right)^{\frac{1}{2}} \left( \intsptiT{(\uet)^{\alpha +3}} \right)^{\frac{1}{2}}  
		\le C
	\end{align}
	for all $\e \in (0,1)$. Combining \eqref{lem:starkeKonvergenz-L4-h2} and \eqref{lem:starkeKonvergenz-L4-h3}, we conclude
	\begin{align}\label{lem:starkeKonvergenz-L4-h4}
		\underset{\e \in (0,1)}{\sup}	\Norm{((\uet)^{\frac{\alpha+3}{2}})}_{H^{2}({\Ot})} < C \,.
	\end{align}
	By Simon's theorem, c.f. \cite{Simon1986}, using \eqref{lem:starkeKonvergenz-L4-h1}, \eqref{lem:starkeKonvergenz-L4-h4}, and the spaces 
	$H^{2}(\cO) \subset H^{1}(\cO) \subset L^{1}(\cO)$, 
	we get
	\begin{align}\label{starkeKonvergenz-uex}
		((\uet)^{\frac{\alpha+3}{2}})_{x} \rightarrow (\ut^{\frac{\alpha+3}{2}})_{x}
	\end{align}
	strongly in $L^{2}(\cO_{T})$.
	On the set $[\ut > 0]$ there exists a subsequence with
	\begin{align}\label{punktweiseKonvergenz-uex}
		((\uet)^{\frac{\alpha + 3}{4}})_{x} \rightarrow (\ut^{\frac{\alpha + 3}{4}})_{x}
	\end{align}
	pointwise almost surely for $\e \rightarrow 0$.
	The sequence $((\uet)^{\frac{\alpha + 3}{4}})_{x}$ is uniformly bounded in $\lpr{4}{\Ot}$, cf. \eqref{est:ApproxSolutions}, which in turn implies uniform integrability in $\lpr{4-\delta}{\Ot}$, $\delta \in (0,1)$. The result now follows with Vitali's theorem. 
	\proofend
	
	\begin{kor1}\label{kor:schwachekonvergenz-uxet^3}
		We have $\ti{\mathbb{P}}$-almost surely 
		\begin{align}\label{kor:schwachekonvergenz-uxet^3h1}
			(\uext)^{3} &\rightharpoonup \uxt^{3}
		\end{align}
		weakly in $L^{\frac{4}{3}}([\ut>0])$
		and for $\phi \in H^{3}_{\text{per}}(\cO)$  
		\begin{align}\label{kor:schwachekonvergenz-uxet^3h2}
			\underset{\e\rightarrow 0}{\lim}	\int_{[\ut=0]}{(\uext)^{3}\phi_{x}} \dx = 0 .
		\end{align}
	\end{kor1}
	\proof 
	The identity (\ref{eq:Darstellung-uxe^3}), the uniform convergence of $\uet$, and Lemma \ref{lem:starkeKonvergenz-L4-} show 
	$ 
	(\uext)^{3} \rightarrow 	\uxt^{3}
	$
	strongly in $L^{\frac{4}{3}-}([\ut>0])$. By the uniqueness of weak limits we find $\zeta = (\uext)^{3}$ on $[\ut>0]$ in Lemma \ref{lem:schwacheKonvergenz-uex^3} which is \eqref{kor:schwachekonvergenz-uxet^3h1}. The second statement \eqref{kor:schwachekonvergenz-uxet^3h2} follows with Hölder's inequality:
	\begin{align}\notag
		&\int \int_{[\ut=0]}\abs{(\uext)^{3}\phi_{x}}\dx \ds =	\int \int_{[\ut=0]}{(\uet)^{\frac{3(\alpha -1)}{4}}\abs{\uext}^{3}(\uet)^{\frac{-3(\alpha -1)}{4}} }\abs{\phi_{x}} \dx \ds \\ \notag
		&\le 	
		\left( \int \int_{[\ut=0]} (\uet)^{\alpha-1}(\uext)^{4}	\dx \ds	\right)	^{3/4}
		\left( \int \int_{[\ut=0]} (\uet)^{-3\alpha+3}\phi_{x}^{4}	\dx \ds	\right)	^{1/4} 	\rightarrow	0 \, .
	\end{align}
	\proofend 
	
	The convergence of $\uet(\uext)^{2}$ can be shown in a similar way. 
	\begin{lem1}\label{lem:schwacheKonvergenzTerm2}
		For a subsequence of $\uet$ we have $\ti{\mathbb{P}}$-almost surely
		\begin{align}
			\uet(\uext)^{2} \rightharpoonup 	\ut(\uxt)^{2} 
		\end{align} 
		weakly in $\lpr{2}{\cO_{T}}$. 
	\end{lem1} 	
	\proof
	By means of the identity 
	\begin{align}\label{lem:schwacheKonvergenzTerm2h1}
		\uet(\uext)^{2} = \left(\frac{4}{\alpha+3}\right)^{2}(\uet)^{\frac{-\alpha+3}{2}}((\uet)^{\frac{\alpha +3}{4}})^{2}_{x} \, ,
	\end{align}
	we conclude as in Lemma \ref{lem:schwacheKonvergenz-uex^3} that $\uet(\uext)^{2}$ is uniformly bounded in $\lpr{2}{\Ot}$ and thus  admits a subsequence such that $\uet(\uext)^{2} \rightharpoonup  \gamma$ $\ti{\mathbb{P}}$-almost surely in $\lpr{2}{\cO_{T}}$.
	Using $-\alpha > 0$, the strong convergence (\ref{starkeKonvergenz-uex}) in Lemma \ref{lem:starkeKonvergenz-L4-}, as well as  
	\begin{align}\label{lem:schwacheKonvergenzTerm2h2}
		\uet (\uext)^{2} = \left(\frac{2}{\alpha +3}\right)^{2} ((\uet)^{\frac{\alpha+3}{2}})^{2}_{x} (\uet)^{-\alpha},
	\end{align}
	we find a subsequence of $\uet (\uext)^{2}$ that converges pointwise to $\ut (\uxt)^{2}$ $\ti{\mathbb{P}}$-almost surely. The uniform bound of \eqref{lem:schwacheKonvergenzTerm2h2} and Vitali's theorem then give 
	\begin{align*}
		\uet (\uext)^{2} \rightarrow \ut (\uxt)^{2}
	\end{align*}
	strongly in $\lpr{2-}{\cO_{T}}$ $\ti{\mathbb{P}}$-almost surely.
	The same arguments as in Corollary \ref{kor:schwachekonvergenz-uxet^3} show $\gamma = \ut (\uxt)^{2}$.
	\proofend 
	
	The next result follows immediately from the strong convergence of $(\uet)^2$ and the weak convergence of $\uext$ in $\lpr{2}{\Ot}$ respectively.
	\begin{lem1}\label{lem:schwacheKonvergenzTerm3}
		For all $\phi \in H^{3}_{\text{per}}(\Ot)$ 
		\begin{align}\label{lem:schwacheKonvergenzTerm3h0}
			\intsptiT{(\uet)^{2}\uet_{x}\phi_{xxx}} \rightarrow \intsptiT{(\ut)^{2}\ut_{x}\phi_{xxx}}
		\end{align}
		holds $\ti{\mathbb{P}}$-almost surely.
	\end{lem1}
	
	Finally, we show that the term which contains the effective interface potential vanishes in the limit.
	\begin{lem1}\label{lem:PotentialTermVanishes}
		For $\e\rightarrow 0$ and all $\phi \in H^{3}_{\text{per}}(\Ot)$ we have
		\begin{align}
			\Ew{\abs{ \intsptiT{\e(\uet)^{2} \uext \W''(\uet) \phi_{x} }  }  } \rightarrow 0.
		\end{align}
	\end{lem1}
	\proof
	Using the weighted version of Young's inequality, we find for $\eta > 0$ and $p \in (2,\infty)$ with $\W(x) = x^{-p}$
	\begin{align}\notag
		\abs{\intsptiT{\e(\uet)^{2}\uext (\uet)^{-p-2} \phi_{x}  }}
		&=\abs{ \frac{1}{p-1} \intsptiT{\e(\uet)^{\frac{-p}{2}} (\uet)^{\frac{-p+2}{2}} \phi_{xx}}} 
		\\ \notag
		&\le \frac{C}{4}\intsptiT{\e^{2} \e^{\eta -1} (\uet)^{-p} \phi^{2}_{xx} } 
		\\\notag 
		&\quad +C  \intsptiT{\e^{1-\eta}	(\uet)^{-p+2}} 
		\\ \notag
		&:= \operatorname{I} + \operatorname{II}.
	\end{align}
	Due to the boundedness of $\phi_{xx}$ and \eqref{est:ApproxSolutions}, we have for $ \operatorname{I}$
	\begin{align}\notag
		C \e^{\eta +1} \Ew{\intsptiT{\W (\uet) \phi_{xx}^{2}}}
		\le
		\e^{\eta +1} C \,  \Ew{\sup_{t\in [0,T]}\intsp{\W(\uet(t))}} 
		\le 
		\e^{\eta +1} C \rightarrow 0 
	\end{align}
	$\mathbb{\ti{P}}$-almost surely.
	For $\operatorname{II}$, we argue with $\delta > 0$ as follows: 
	\begin{align}\notag
		& \intsptiT{\e^{1-\eta} (\uet)^{-p}(\uet)^{2}	} 
		\\ \notag
		=&\int \int_{[\ut<\e^{\frac{\eta}{2}+ \delta }]}{\e^{1-\eta} (\uet)^{-p}(\uet)^{2}	} \dx \ds
		+ \int \int_{[\ut \ge \e^{\frac{\eta}{2}+ \delta }]} \e^{1-\eta} (\uet)^{-p+2} \dx \ds
		\\ \notag
		\le &\e^{1-\eta+\eta + 2\delta} \intsptiT{\W(\uet)} 
		+ \int \int_{[\ut \ge \e^{\frac{\eta}{2}+ \delta }]}{\e^{1-\eta+(-p+2)(\frac{\eta}{2}+ \delta)}	}\dx \ds \, .
	\end{align}
	Hence, using \eqref{est:ApproxSolutions} once more, we have for $\eta$ and $\delta$ chosen appropriately
	\begin{align*}
		\Ew{\operatorname{II}}
		\le & \, \e^{1 +2\delta} C 
		+  \e^{1-\eta- \frac{p\eta}{2}+\eta + \delta(- p+ 2)} C
		\rightarrow 0
	\end{align*}
	for $\e \rightarrow 0$.
	\proofend
\subsection{Identification of the Stochastic Integral}
	\label{sec:IdentificationStochInt}
	
	For  $\phi \in H^{3}_{\text{per}}(\cO)$ arbitrary, but fixed, we consider the processes 
	$\M_{\e,\phi}: \Omega \times [0,T] \rightarrow \mathbb{R} $ defined by
	\begin{align}\label{def:Meps}
		\M_{\e,\phi}(t) &:= 
		\intsp{(\ue(t)-u_{0}^{\e})\phi} \notag
		- \intsptit{(\ue_{x})^3\phi_{x}} \notag
		- 3 \intsptit{\ue (\ue_{x})^{2}\phi_{xx}} 
		\\ \notag
		&- \intsptit{(\ue)^{2}\ue_{x}\phi_{xxx}}	
		+ \intsptit{(\ue)^{2} \uex \e \W''(\ue) \phi_{x} }   
		\\
		&+ \intsptit{\frac{1}{2} \sum_{k \in \Z} \lal^{2}\gl(\gl\ue)_{x} \phi_{x} } \, . 
	\end{align}
Note that the right-hand side of \eqref{def:Meps} coincides with the deterministic terms in \eqref{eq:SolutionApproxProblems1} for the choice $\phi\in H^3_{per}(\ort)$ which follows easily by integration by parts. In particular, the last term in \eqref{def:Meps} is identical with $C_{Strat}\int_0^t\int_\ort \ue_x\phi_x dx ds$, cf. \eqref{eq:stratito1}, \eqref{eq:appendix12}, and \eqref{eq:appendix13}.
On the other hand, we have  
	\begin{align}
		\M_{\e,\phi}(t) 	&=  \sum_{k \in \Z} \int_{0}^{t} \int_{\cO} \lal (\gl\ue)_{x}\phi \dx \dbl 
	\end{align}
	for $t \in [0,T]$, i.e. $\M_{\e,\phi}$ is a continuous, square integrable $\mathcal{F}_{t}^{\e}$-martingale.
	We will need the following results which have been shown in \cite{FischerGruen}, Lemma 5.10 and Lemma 5.12:  
	\begin{align}\label{eq:QuadrVarMeps}
		\qv{\M_{\e,\phi}} 
		= \int_{0}^{\cdot}{ \sum_{k \in \Z} \lal^{2} 
			\left(
			\intsp{(\ue \gl)_{x}\phi}
			\right)^{2}	} \ds ,
	\end{align}
	\begin{align}\label{est:QuadrVarMeps}
		\qv{\M_{\e,\phi}} \le C \Norm{\phi}^{2}_{H^{1}_{\text{per}}} \inttimeT{\Norm{\ue}^{2}_{L^{2}(\cO)}},
	\end{align}
	and for $k \in\N$ 
	\begin{align}\label{QuadrCovarMbeta}
		\qv{\M_{\e,\phi}, \ble} = 
		\lal \int^{\cdot}_{0}\int_{\cO}{ (\ue \gl )_{x} \phi} \dx \ds .
	\end{align}
	With these results at hand we can establish 
	\begin{kor1}\label{kor:Martingal1}
		Let $k \in\N$ and $\e \in (0,1)$. The processes
		\begin{align}
			\M_{\e,\phi}^{2} - \int_{0}^{\cdot}{ \sum_{k \in \Z} \lal^{2} 
				\left(
				\intsp{(\ue \gl)_{x}\phi}
				\right)^{2}	} \ds 
		\end{align}
		and 
		\begin{align}
			\M_{\e,\phi}\ble -	\lal \int_{0}^{\cdot}\intsp{ (\ue \gl )_{x} \phi} \ds
		\end{align}	
		are continuous $\mathcal{F}^{\e}_{t}$-martingales. 
	\end{kor1}
	By the equality of laws stated in Theorem \ref{sat:JakubowskiAnwendung}, we also get the analog statements for  
	\begin{align}
		\ti{\M}_{\e,\phi}(t) &:= 
		\intsp{(\uet(t)-\ut_{0}^{\e})\phi} \notag
		- \intsptit{(\uet_{x})^3\phi_{x}} \notag\\ \notag
		&- 3 \intsptit{\uet (\uet_{x})^{2}\phi_{xx}} 
		- \intsptit{(\uet)^{2}\uet_{x}\phi_{xxx}}	\\ 
		&+ \intsptit{(\uet)^{2} \uext \e \W''(\ue) \phi_{x} }  
		+ \intsptit{\frac{1}{2} \sum_{k \in \Z} \lal^{2}\gl(\gl\uet)_{x} \phi_{x} } \, .
	\end{align}

	\begin{lem1}\label{sat:martingaleigenschaftTilde} 
		For $k \in\N$ and $\e \in (0,1)$
		\begin{align}\label{sat:martingaleigenschaftTildehh1}
			&\ti{\M}_{\e,\phi}
			\\ 
			\label{sat:martingaleigenschaftTildehh2}
			&\ti{\M}_{\e,\phi}^{2} - \int_{0}^{\cdot}{ \sum_{k \in \Z} \lal^{2} 
				\left(
				\intsp{(\uet \gl)_{x}\phi}
				\right)^{2}	} \ds  
			\\ \label{sat:martingaleigenschaftTildehh3}
			&\ti{\M}_{\e,\phi}\ti{\ble} -	\lal \int_{0}^{\cdot}\int_{\cO} (\uet \gl )_{x} \phi \dx \ds 
		\end{align}
		are continuous $\ti{\mathcal{F}}_{t}^{\e}$-martingales. Moreover, on $[0,T]$ we have
		\begin{align}\label{sat:martingaleigenschaftTildeh1}
			&\qv{\ti{\M}_{\e,\phi}}_{t} =  \int_{0}^{t}{ \sum_{k \in \Z} \lal^{2} 
				\left(
				\intsp{(\uet \gl)_{x}\phi}
				\right)^{2}	} 
			\ds 
			\\\label{sat:martingaleigenschaftTildeh2}
			&\qv{\ti{\M}_{\e,\phi},\ti{\ble}}_{t} = 	\lal \intsptit{ (\uet \gl )_{x} \phi}.
		\end{align}
	\end{lem1}
	
	The next step is to show that the martingale property is preserved in the limit. We show that for $\phi \in H^{3}_{\text{per}}(\cO)$
	\begin{align}\label{def:M0t}
		\ti{\M}_{0,\phi}(t) &:=
		\intsp{(\ut(t)-\unt)\phi} \notag
		- \int \int_{[\ut>0]}(\ut_{x})^3\phi_{x}  \dx \ds \notag
		- 3 \intsptit{\ut (\ut_{x})^{2}\phi_{xx}} \\ 
		&- \intsptit{(\ut)^{2}\ut_{x}\phi_{xxx}}	
		+ \intsptit{\frac{1}{2} \sum_{k \in \Z} \lal^{2}\gl(\gl\ut)_{x} \phi_{x} }
	\end{align}
	has the martingale property.
	\begin{lem1}\label{sat:M0Martingal}
		For $s,t \in [0,T]$ with $s \le t$ and
		for all continuous functions $\Psi \colon C^{\ti{\gamma} ,\ti{\sigma}}([0,s]\times \bar{\cO)}\times C([0,s];L^{2}(\cO)) \rightarrow [0,1]$, we have 
		\begin{align}\label{sat:M0Martingalh1}
			\Ew{\Psi(r_{s}\ut, r_{s} \ti{W}) \left(\ti{\M}_{0,\phi}(t)-\ti{\M}_{0,\phi}(s)\right)} = 0. 
		\end{align}
	\end{lem1}  
	\proof 
	We treat the terms inside the expectation in (\ref{sat:M0Martingalh1}) one by one. The continuity of $\Psi$ as well as the convergence of $\uet$ to $\ut$ and of $\ti{W}^{\e}$ to $\ti{W}$ in $C(0,T;L^{2}(\cO))$ show
	\begin{align}\label{KonvergenzPsi}
		\underset{\e \rightarrow 0}{\lim} \Psi(r_{s}\uet, r_{s} \ti{W}^{\e}) =  \Psi(r_{s}\ut, r_{s} \ti{W})
	\end{align}
	$\ti{\mathbb{P}}$-almost surely on $[0,1]$. To see the convergence of the expected values, we aim to utilize Vitali's theorem; therefore, since $\Psi$ is bounded, it suffices to show uniform boundedness of moments of the integral-terms in (\ref{sat:M0Martingalh1}) and use the convergence results already established.
	
	By the strong convergence of $\uet$ in $ C^{\ti{\gamma} ,\ti{\sigma}}(\Otb)$, cf. Theorem \ref{sat:JakubowskiAnwendung}, we have $\ti{\mathbb{P}}$-almost surely 
	\begin{align}\label{Konvergenzue}
		\underset{\e \rightarrow 0}{\lim}	\intsp{(\uet(t)-\uet(s))\phi} =  \intsp{(\ut(t)-\ut(s))\phi} \, .
	\end{align}    
	The $\alpha$-entropy-energy estimate \eqref{est:ApproxSolutions} gives
	$	
	\uet \in \lpr{2q}{\ti{\Omega};\lpr{\infty}{\cO_{T}}}
	$
	for an arbitrary $q > 1$ and thus, the uniform boundedness of a $q$-th absolute moment.
	
	Weak convergence of $(\uext)^{3}$ $\mathbb{\ti{P}}$-almost surely has been established in Corollary \ref{kor:schwachekonvergenz-uxet^3}.
	Using Hölder's inequality and \eqref{est:ApproxSolutions}, we have for $q>1$ 
	\begin{align*}\notag
		&\Ew{\abs{\intsptist{(\uext)^{3}\phi_{x} } }^{q} }
		\\ \notag	
		\le
		C \, &\Ew{\abs{\intsptist{ ((\uet)^{\frac{\alpha +3}{4}})^{4}_{x}} }^{q} }^{\frac{3}{4}}	
		\, \Ew{\abs{\intsptist{(\uet)^{-3(\alpha-1)}  \phi_{x}^{4} } }^{q} } ^{\frac{1}{4}} 
		\le C.
	\end{align*} 
	
	The identity
	\begin{align*}
		(\uext)^{2}\uet = \left(\frac{4}{\alpha +3} \right)^{2} ((\uet)^{\frac{\alpha+3}{4}})_{x}^{2} (\uet)^{\frac{-\alpha+3}{2}}
	\end{align*}
	yields
	\begin{align*}
		&\Ew{\abs{\intsptist{(\uext)^{2}\uet\phi_{xx}} }^{q}  }
		\\ 
		\le 
		C \, &\Ew{\abs{\intsptist{ ((\uet)^{\frac{\alpha +3}{4}})^{4}_{x}} }^{q} }^{\frac{1}{2}}	
		\, \Ew{\abs{\intsptist{(\uet)^{-\alpha +3}\phi_{xx}^{2} } }^{q} } ^{\frac{1}{2}}
		\le 
		C 
	\end{align*}
	which we combine with Lemma \ref{lem:schwacheKonvergenzTerm2}.
	
	By means of Cauchy-Schwarz' inequality 
	\begin{align*}
		\abs{\intsptist{(\uet)^{2}\uet_{x}\phi_{xxx}}}  \le C
		\left(\underset{\Otb}{\sup} \, \uet \right)^{2}  \left(\underset{t \in [0,T]}{\sup} \intsp{\abs{\uext}^{2}}\right)^{\frac{1}{2}}
		\left( \intsp{\phi_{xxx}^{2}}\right)^{\frac{1}{2}}
	\end{align*}
	holds, which implies the boundedness of higher moments. Lemma \ref{lem:schwacheKonvergenzTerm3} states the needed convergence $\ti{\mathbb{P}}$-almost surely.
	
	In Lemma \ref{lem:PotentialTermVanishes} we have seen
	\begin{align}\label{sat:M0Martingalh4}
		\underset{\e \rightarrow 0}{\lim} \,	\Ew{\Psi(r_{s}\uet, r_{s} \ti{W}^{\e}) \left( \intsptiT{\e(\uet)^{2} \uex \W''(\uet) \phi_{x} } \right) } = 0 \, .
	\end{align}
	For the Stratonovich correction term, we have due to the convergence of $\uet$ $\ti{\mathbb{P}}$-almost surely and the weak convergence of $\uext$ on $\lpr{2}{0,T;\lpr{2}{\cO}}$, cf. Theorem \ref{sat:JakubowskiAnwendung},
	\begin{align*}\notag
		\intsptist{\frac{1}{2} \sum_{k \in \Z} \lal^{2}\gl(\gl\uet)_{x} \phi_{x} }\rightarrow
		\intsptist{\frac{1}{2} \sum_{k \in \Z} \lal^{2}\gl(\gl\ut)_{x} \phi_{x} }
	\end{align*}
	$\ti{\mathbb{P}}$-almost surely. Furthermore, by boundedness of the $\gl$ and \eqref{Lambda-Wachstumsbedingung},   
	$$
	\abs{\intsptist{\frac{1}{2} \sum_{k \in \Z} \lal^{2}\gl(\gl\uet)_{x} \phi_{x} }}\le	C \abs{\underset{\Otb}{\sup}\, \uet}.
$$ 
	
	From Lemma \ref{sat:martingaleigenschaftTilde} we know that $(\ti{\M}_{\e,\phi})_{\e \in (0,1)}$ are martingales. Thus, 
	the convergence results above yield
	\begin{align*}
		\Ew{\Psi(r_{s}\ut, r_{s} \Wt) \left( \ti{\M}_{0,\phi}(t)-\ti{\M}_{0,\phi}(s)\right)} = 0.
	\end{align*}
	\proofend
	
	Dynkin's lemma in combination with Lemma \ref{sat:M0Martingal} implies the martingale property, cf. for example \cite{HofmanovaSeidler2012}.
	\begin{kor1}\label{kor:M0Martingal}
		$\ti{\M}_{0,\phi}$ is a continuous $\ti{\F_{t}}$-martingale.
	\end{kor1}
	
	By similar arguments as before, cf. also \cite{FischerGruen} Lemmas 5.14 and 5.15, we can show that for 
	$0 \le s \le t \le T$ and $\Psi$ as in Lemma \ref{sat:M0Martingal}
	\begin{align}\label{sat:Martingalea1}
		\Ew{ \Psi(r_{s}\ut,r_{s}\ti{W})  \left( \ti{\M}_{0,\phi}^{2}(t)- \ti{\M}_{0,\phi}^{2}(s) - \int_{s}^{t}{ \sum_{k \in \Z} \lal^{2} 
				\left(
				\intsp{(\ut \gl)_{x}\phi}
				\right)^{2}	} \ds \right)} = 0
	\end{align}
	and 
	\begin{align}\label{sat:Martingalea2}
		\Ew{ \Psi(r_{s}\ut,r_{s}\ti{W}) \left( (\ti{\M}_{0,\phi}\blt)(t) - (\ti{\M}_{0,\phi}\blt)(s) -	\lal \intsptist{ (\ut \gl )_{x} \phi}\right)  } = 0
	\end{align}
	holds.

Following the argumentation of Lemma 5.16 in \cite{FischerGruen}, the identification of the stochastic term is achieved. 
	\begin{lem1}\label{sat:IdentM0}
		It holds $\mathbb{\ti{P}}$-almost surely
		\begin{align}\label{sat:IdentM0-h3}
			\ti{\mathcal{M}}_{0,\phi} =  \sum_{k \in \Z}\int_{0}^{\cdot}\intsp{ \lal (\ut \gl)_{x}\phi} \dblt \, .
		\end{align}
		
	\end{lem1}

        \subsection{Proof of the Main Results}
	Finally, we provide the proofs for the existence of zero-contact angle martingale solutions.
	\begin{myproof}{ of Theorem \ref{theo:MainResult}}
	From Theorem \ref{sat:JakubowskiAnwendung} the existence of the stochastic basis, the $Q$-Wiener process $\ti{W}$, as well as sequences $(\uet)_{\e \in (0,1)}$, $(\uent)_{\e \in (0,1)}$, and random variables $\ut$ and $\unt$ follows. Moreover, for every $\e \in (0,1)$ $\ue$ and $\uet$ as well as $\uen$ and $\uent$ have the same laws, respectively, and
	for $\e \rightarrow 0$,
	$
	\uet \rightarrow \ut
	$
	in $C^{\ti{\gamma} ,\ti{\sigma}}(\Otb)$ and
	$
	\uent \rightarrow \unt
	$
	in $H^{1}_{\text{per}}(\cO)$ holds  $\ti{\mathbb{P}}$-almost surely for a subsequence.
	Owing to the uniform convergence $\uet \rightarrow \ut$, we see in particular that $\ut$ is nonnegative $\ti{\mathbb{P}}$-almost surely.
	The weak formulation \eqref{eq:MainResultWeakForm} is satisfied due to Lemma \ref{sat:IdentM0}.
	Since $\un + \e^{\theta} \rightarrow \un$ pointwise, by Corollary 13.19 in \cite{Klenke2014eng} we get 
	\begin{align*}
		\mathbb{P}^{\e} \circ (\uen)^{-1} = \mathbb{P} \circ (\un + \e^{\theta})^{-1} \rightarrow \mathbb{P} \circ \un^{-1} = \Lambda^{0}
	\end{align*}
	weakly for $\e \rightarrow 0$.
	Likewise, from the pointwise convergence $\uent \rightarrow \unt$ $\ti{\mathbb{P}}$-almost surely, the weak convergence of the laws follows, i.e. 
	\begin{align*}
		\ti{\mathbb{P}} \circ (\uent)^{-1} \rightarrow \ti{\mathbb{P}} \circ \unt^{-1}.
	\end{align*}
	From Theorem \ref{sat:JakubowskiAnwendung} we know that
	\begin{align*}
		\mathbb{P}^{\e} \circ (\uen)^{-1} = \ti{\mathbb{P}} \circ (\uent)^{-1}
	\end{align*}
	for all $\e \in (0,1)$. Thus, by the uniqueness of limits w.r.t. weak convergence on Polish spaces, we get
	\begin{align*}
		\Lambda^{0} = \ti{\mathbb{P}} \circ \unt^{-1} \, .
	\end{align*}
	Additionally, $\unt$ is $\ti{\mathbb{P}}$-almost surely nonnegative.
	
	At last, estimate \eqref{est:MainResult} follows from Fatou's lemma and \eqref{est:ApproxSolutions}. We have
	\begin{align}
		\notag
		\Ew{\underset{\e \rightarrow 0}{\liminf} \,\sup_{t\in [0,T]} \left (\intsp{\frac{1}{2} \abs{\uet_{x}}^{2}}  \right)^{q}}
		\le
		\notag
		\underset{\e \rightarrow 0}{\liminf} \, \Ew{\sup_{t\in [0,T]} \left (\intsp{\frac{1}{2} \abs{\uet_{x}}^{2}}  \right)^{q}} 
		\le  C(\ut_{0},q,T).
	\end{align}
	
	Hence, by the lower semi-continuity of the $\lpr{\infty}{0,T;H^{1}(\cO)}$ norm w.r.t. the convergence in the distributional sense and Remark 8.3 from \cite{Alt2016eng}, we deduce 
	\begin{align*}
		\Ew{  \sup_{t\in [0,T]}\left(\frac{1}{2} \intsp{\abs{\uxt}^{2} } \right)^{q} } \le  C(\ut_{0},q,T) \, . 
	\end{align*}
	Again, with Fatou's lemma and the lower semi-continuity of the $L^{p}(\Ot)$-norm w.r.t. to weak convergence of $((\uet)^{\frac{\alpha +3 }{4}})_{x}$ and $((\uet)^{\frac{\alpha +3 }{2}})_{xx}$, we conclude 
	\begin{align*}
		\Ew{\left(\intsptiT{((\ut)^{\frac{\alpha + 3}{4}})_{x}^{4}}\right)^{q}}
		\le C(\ut_{0},q,T)
	\end{align*}
	and 
	\begin{align*}
		\Ew{\left(\intsptiT{((\ut)^{\frac{\alpha+3}{2}})^{2}_{xx}}\right)^{q}} 
		\le C(\ut_{0},q,T),
	\end{align*}
	where we used the $\alpha$-entropy-energy estimate respectively.
	Thus, we have shown the estimate (\ref{est:MainResult}) and completed the proof of Theorem \ref{theo:MainResult}. 
	\end{myproof}
Finally, we show the zero-contact-angle property at touch-down points of the solution.
\begin{myproof}{ of Corollary~\ref{cor:contactangle}}
From the a priori estimate \eqref{est:MainResult} combined with (H2), we infer that $\ut^{\tfrac{\alpha+3}{4}}(\cdot,\cdot,\omega)$ is element of $L^4((0,T);W^{1,4}(\ort))$ $\tprobab$-almost surely. Using the nonnegativity of $\ut$, the assumption on $\alpha$,  and the fact that $\ut^{\tfrac{\alpha+3}{4}}(\cdot,t,\omega)\in W^{1,4}(\ort)$ for almost all $t\in [0,T]$, the claim follows from the estimate
$$0\leq \ut(x,t_0, \omega)\leq C(\omega,t_0)|x-x_0|^{\tfrac{3}{\alpha+3}} \qquad \text{ for }x\in\ort$$ which  is a consequence of Sobolev's embedding theorem.
\end{myproof}

{\bf Concluding remarks.} In this paper which is partially based on the master thesis of the second author \cite{KleinMaster}, we have presented a rather elementary proof for the existence of zero-contact angle solutions to the stochastic thin-film equation \eqref{eq:STFES} for a quadratic mobility $m(\cdot)$ in one space dimension. The strategy has been to derive new regularity results first for approximate solutions which differ from \eqref{eq:STFES} by a potential that enhances spreading and that this way entails strict positivity almost surely. We expect that this method can be slightly modified to establish corresponding results in the spatially two-dimensional case, too, this time starting from the existence result in \cite{MetzgerGruen2021}. Moreover, the new integral estimate \eqref{est:MainResult} may serve as a starting point to establish results on the qualitative behaviour of solutions -- like finite speed of propagation or (non)-occurrence of waiting time phenomena. It remains, however, an open problem to which extent this approach may be applied to more general mobilities $m(\cdot)$. 

	\appendix
	
\section{Stratonovich Correction}
	\label{StratCorrection}
	We will briefly discuss how to derive equation \eqref{eq:stratito1} from equation \eqref{eq:StfeStrato}.
	We skip the index $\e$ in this section. The Stratonovich correction term with respect to
	\begin{align}
		\sum_{k \in \Z} (\lal \gl u)_{x} \circ \dbl(t)
	\end{align} 
	reads 
	\begin{align}
		\mathcal{C}_{S} =  \frac{1}{2} \sum_{ k \in \Z} \lal^{2} (\gl  (u \gl )_{x} )_{x} \, ,
	\end{align}
	see, e.g., \cite{GessGnann2020}.
	We list some frequently used identities which follow by direct computation from \eqref{def:ONB} and \eqref{def:Eigenvalues}:
	\begin{align}\label{eq:ONB-1}
		\gl \gl ' 
		&=
		\frac{2\pi k}{L} \gl \gml	 
		\\	\label{eq:ONB-2}
		\gml  \gml' 
		&=
		-\frac{2\pi k}{L} \gml \gl	 
		\\	\label{eq:ONB-3}
		\gl^{2} + \gml^{2} 
		&=
		\frac{2}{L}	  
	\end{align}
	and consequently
	\begin{align}
		\label{eq:ONBRelations1}
		\sum_{k \in \Z}\lal^{2}     \gl^{2}  
		&=
		\frac{\lan^{2}}{L}
		+
		\sum_{k =1}^{\infty} \frac{2\lal^{2}}{L}
		\\	\label{eq:ONBRelations2}
		\sum_{k \in \Z} \lal^{2}   \glx^{2}  
		&=
		8 \pi^2\sum_{k =1}^{\infty} \frac{\lal^{2}k^{2}}{L^{3}}
		\\ \label{eq:ONBRelations3}
		\sum_{k \in \Z} \lal^{2}    \glx\gl  
		&=
		0
		\\	\label{eq:ONBRelations4}
		\sum_{k \in \Z}  \lal^{2}   \glxx^{2}  
		&=
		32 \pi^{4}\sum_{k =1}^{\infty} \frac{\lal^{2}k^{4}}{L^{5}}
		\\	\label{eq:ONBRelations5}
		\sum_{k \in \Z}   \lal^{2}      \glx\glxx  
		&=
		0
		\\	\label{eq:ONBRelations6}
		\sum_{k \in \Z}  \lal^{2}      \glxx \gl  
		&=
		- 8\pi^2
		\sum_{k =1}^{\infty} \frac{\lal^{2}k^{2}}{L^{3}} \, .
	\end{align}

	We have
	\begin{align} \notag
		(\gl (\gl u)_{x})_{x} 
		&= \glx^{2}u + \gl \glxx u + 3\gl \glx \ux + \gl^{2} \uxx \, .
	\end{align}
	Using the relations above, we therefore get
	\begin{equation} \label{eq:appendix12} 
		\sum_{k \in \Z} \lal^{2}(\gl (\gl u)_{x})_{x} 
		=  
		\left(
		\frac{\lan^{2}}{L}
		+
		\sum_{k = 1}^{\infty} \frac{2\lal^{2}}{L}
		\right)
		u_{xx} \, .
	\end{equation}
	Hence, 
	\begin{align}\label{eq:appendix13}
		\mathcal{C}_{S} = \Cstr u_{xx} \, .
	\end{align} 
	The stochastic thin-film equation with Stratonovich noise can then  be written as
	
	\begin{equation}
		\dif u=(-(u^2u_{xxx})_x+C_{Strat} u_{xx}) \dt +(u \dif W)_x 
	\end{equation}
	or equivalently	as
	\begin{align} 
			\dif u &=  - ( u^{2}(u_{xx} - \Cstr(1- u^{-1}) )_{x} )_{x} \dt  
			+ (u\dif W)_x
			\\
			&=
			 -( u^2(u_{xx}-\mathcal S'(u))_x)_x\dif t +(u\dif W)_x
	\end{align}		
 	with $\mathcal{S}(u) = \Cstr (u - \log u)$.
 	
\section{It\^o's Formula}\label{appendix:Itoformula}
	In what follows, we will first state It\^o's formula from \cite{Krylov} that is used in the proof of Lemma \ref{lem:alpha-Entropy-Energy}. Subsequently, we will justify its application in our setting.
	\begin{sat1}\label{theo:Ito/KrlyovFormel}
		Let $V$
		and $H$ be separable Hilbert-spaces such that $V$ is a dense subset of $H$ in the metric of $H$ and assume $\Norm{u}_{H} \leq \Norm{u}_{V}$ for all $u \in V$. 
		Furthermore, consider a complete probability space $(\Omega , \mathcal{F}, \mathbb{P})$ with an increasing complete filtration
		$(\mathcal{F}_{t})_{t\ge0} \subset \mathcal{F}$
		and predictable stochastic processes $v$ and $v^{\ast}$ that satisfy for all $T >0$
		\begin{align}\label{Krylov:Processes}
			\Ew{\int_{0}^{T}\Norm{v_{t}}^{2}_{V} + \Norm{v^{\ast}_{t}}^{2}_{V} \dt } <\infty  
		\end{align}
		and for all $\phi \in V$, $t \in [0,T]$ 
		\begin{align}\label{Krylov:WeakForm}
			\skp{\phi}{v_{t}}_{H} = \skp{\phi}{v_{0}}_{H} 
			+
			\int_{0}^{t}\skp{\phi}{v_{s}^{\ast}}_{V} \ds
			+ 
			\sum_{k} \inttimetra{\skp{\phi}{\ssl}_{H}} \, ,
		\end{align}
		where $\ssl$ are predictable $H$-valued processes with 
		\begin{align}\label{Krylov:StochInt}
			\sum_{k\in \Z} \Ew{\int_{0}^{T} \Norm{\ssl}^{2}_{H} \dt } < \infty 
		\end{align} 
		for every $T > 0$. Let $\Phi$ be a real-valued function on $H$ with the following properties holding $\mathbb{P}$-almost surely:
		\begin{enumerate}
			\item[i)]
			For every $h,\xi \in H$ the function $\Phi (h +t\xi)$  is twice continuously differentiable as a function of $t$ and the mappings
			\begin{align}\label{Krylov:FrechDer}
				\Phi_{(\xi)}(h) := \left. \abl{t} \Phi(h+t\xi)\right\vert_{t=0}, \quad \Phi_{(\xi)(\xi)}(h) := \left. \ablsq{t} \Phi(h+t\xi)\right\vert_{t=0} 
			\end{align}
			are continuous as functions of $(h,\xi) \in H \times H$.
			\item[ii)]
			For any $R\in (0,\infty)$ there is a constant $K(R)$ such that for all $h, \xi \in H$ with $\Norm{h}_{H} \leq R$ 
			\begin{align}\label{Krylov:LinearQuadraticEst}
				|\Phi_{(\xi)}(h)| \leq K(R) \Norm{\xi}_{H}, \quad 
				|\Phi_{(\xi)(\xi)}(h)| \leq K(R) \Normsq{\xi}_{H}  
			\end{align}
			hold. 
			\item[iii)]
			For $ h \in V $ and a constant $K_{1}$, we have $ \Phi_{(\cdot)}(h) \in V $ and 
			\begin{align}\label{Krylov:RieszProjBound}
				\Norm{\Phi_{(\cdot)}(h)}_{V} \leq K_{1}(1+\Norm{h}_{V}) \,, 
			\end{align}
			where $ \Phi_{(\cdot)}(h)$ is the unique element in $H$   such that the continuous linear functional   $\Phi_{(\xi)}(h)$ can be written as
			\begin{align}\label{Krylov:RieszProjDef}
				\Phi_{(\xi)}(h) = \skp{\Phi_{(\cdot)}(h)}{\xi}_{H}, \quad \Norm{\Phi_{(\cdot)}(h)}_{H} \le K \Norm{h}_{H} \, .
			\end{align}	
			\item[iv)]
			For every $v^{*}\in V$ the function $\skp{\Phi_{(\cdot)}(v)}{v^{*}}_{V}$ is continuous on $V$ in the metric of $V$.
		\end{enumerate}	 
		Then, we have
		\begin{align}\label{Krlyov-Formel}\notag
			\Phi(v_{t}) = \Phi(v_{0})  
			&+ \inttimet{\skp{\Phi_{(\cdot)}(v_{s})}{v_{s}^{*}}_{V}}  \notag
			+ \inttimet{\frac{1}{2} \sum_{k \in \Z}\Phi_{(\ssl)(\ssl)}(v_{s})} \\
			&+ \sum_{k \in \Z}\int_{0}^{t}\Phi_{(\ssl)}(v_{s}) \dif \beta^{k}_{s} \, .
		\end{align}
	\end{sat1}
	Let us first derive weak formulations as in \eqref{Krylov:WeakForm} for the spaces  $V_{1} := H^{2}_{\text{per}}(\cO)$, $V_{2} := H^{1}_{\text{per}}(\cO)$, and $H := \lpr{2}{\cO}$. 
	By Lemma \ref{lem:Minimum-ue}, we see that for $t \in [0,T_{\sigma}]$
	\begin{align}\label{est:ue-positive}
		\ue(\cdot, t) \ge \bar{c}_{\e} \sigma^{\frac{2}{p-2}},
	\end{align}
	where $\bar{c}_{\e} :=\bar{C}_{p} \e^{\frac{1}{p-2}} $.
	For the functions $\us$ as introduced in \eqref{def:ModifiedSolutions}, we get for all $\phi \in H^{1}_{\text{per}}(\cO)$
	\begin{align}\label{eq:ModifiedWeakFormulation}\notag
		\int_{\cO} \us (t)\phi \dx 
		=
		&\int_{\cO}\uen \phi \dx + \int_{0}^{t \wedge T_{\sigma}} \int_{\cO}(-(\us)^{2} \psx - C_{Strat} \usx) \phi_{x} \dx \ds
		\\
		&- 
		\sum_{k \in \Z} \int_{0}^{t \wedge T_{\sigma}} \int_{\cO} \us \gl \phi_{x} \dx \dbl
	\end{align}
	on $[0,\Tmax]$.
	
	By estimate \eqref{lem:UniformEstimate-u^2ph1}, the function $-(\us)^{2} \psx - C_{Strat} \usx$ defines a mapping 
	\begin{align}
		\xi \in \lpr{2}{\Omega \times [0,\Tmax]; (H^{2}_{\text{per}}(\cO))'}
	\end{align} 
	by 
	\begin{align}
		v \mapsto \Ew{\int_{0}^{\Tmax}\int_{\cO} (-(\us)^{2} \psx - \Cstr \usx)_{xx} v \dx \ds} \, .
	\end{align}
	Hence, by Riesz' representation theorem there is $\xiast \in \lpr{2}{\Omega \times [0,\Tmax]; H^{2}_{\text{per}}(\cO)}$ such that
	\begin{align}\label{RieszProjH2}
		\int_{\Omega}
		\int_{0}^{\Tmax}
		{}_{(H^{2}_{\text{per}}(\cO))'}
		\langle 
		\xi, \phi  
		\rangle_{H^{2}_{\text{per}}(\cO)}
		\ds
		\dif\mathbb{P}^{\e}
		=
		\int_{\Omega}
		\int_{0}^{\Tmax}
		\skp{\xiast}{\phi}_{\Hsob{2}{\cO}}
		\ds	
		\dif\mathbb{P}^{\e}
	\end{align} 
	holds for all $\phi \in \lpr{2}{\Omega \times [0,\Tmax]; H^{2}_{\text{per}}(\cO)}$.
	Similarly, we introduce $\fast$ w.r.t. $H^{1}_{\text{per}}(\cO)$, i.e. $\fast$ solves 
	\begin{align}\label{RieszProjH1}
		\int_{\Omega}
		\int_{0}^{\Tmax}
		{}_{(H^{1}_{\text{per}}(\cO))'}
		\langle 
		f, \phi  
		\rangle_{H^{1}_{\text{per}}(\cO)} 
		\ds
		\dif\mathbb{P}^{\e}
		=
		\int_{\Omega}
		\int_{0}^{\Tmax}
		\skp{\fast}{\phi}_{\Hsob{1}{\cO}}
		\ds	
		\dif\mathbb{P}^{\e}
	\end{align}
	for every $\phi \in \lpr{2}{\Omega \times [0,\Tmax]; H^{1}_{\text{per}}(\cO)}$,
	where $f \in \lpr{2}{\Omega \times [0,\Tmax]; (H^{1}_{\text{per}}(\cO))'}$ is  defined via
	\begin{align}
		v \mapsto \Ew{\int_{0}^{\Tmax}\int_{\cO} (-(\us)^{2} \psx - \Cstr \usx) v_{x} \dx \ds} \, . 
	\end{align}
	Riesz' representation theorem also shows
	\begin{align}\notag
		&\Ew{
			\int_{0}^{\Tmax}\left(\Norm{\us}^{2}_{H^{1}(\cO)} 
			+
			\Norm{\fast}^{2}_{H^{1}(\cO)}	\right) \ds 
		} 
		\\ 
		\le
		\Tmax &\Ew{\underset{t \in [0,T]}{\sup}\mathcal{E}_{\e}[\us]} 
		+
		\Ew{\int_{0}^{T}\int_{\cO}(\us^{2}(\ps)_{x})^{2}\dx\ds }
		\le 
		C(\e,T) \, ,
	\end{align}
	due to \eqref{lem:UniformEstimate-u^2ph1}, so \eqref{Krylov:Processes} is satisfied in this case.
	Owing to the regularity of solutions, cf. Definition \ref{def:weak-martingale-solution} iii),  
	\begin{align}
		\Ew{\int_{0}^{\Tmax}
			\left(
			\Norm{\usx}^{2}_{H^{2}(\cO)} 
			+
			\Norm{\xiast}^{2}_{H^{2}(\cO)}	
			\right) \ds } 
		\le
		C(\e,T) 
	\end{align}
	holds as well. 
	
	Let us set $\ssl := (\lal \ue \gl)_{x}$.
	Regarding the predictability of the processes $\xiast$, $\fast$, $\ssl$, and $(\ssl)_{x}$, we refer to Appendix \ref{predictability}.
	Thus, we may rewrite \eqref{eq:ModifiedWeakFormulation} by means of \eqref{RieszProjH2} and obtain
	\begin{align}\label{weakformF1} \notag
		\skp{(\us)_{x}(t)}{\phi}_{\lpr{2}{\cO}} 
		=
		\skp{(\uen)_{x}}{\phi}_{\lpr{2}{\cO}}
		&+
		\int_{0}^{t \wedge T_{\sigma}} \skp{-\xi^{\ast} }{\phi}_{H^{2}(\cO)} \ds
		\\
		&+
		\sum_{ k \in \Z} \int_{0}^{t\wedge T_{\sigma}}  \skp{(\ssl)_{x}}{\phi}_{\lpr{2}{\cO}}  \dbl \, , 
	\end{align}
	where we have multiplied with $\phi_{x}$ for $\phi \in \Hsobper{2}{\cO}$ and integrated by parts. On the other hand, with \eqref{RieszProjH1} we get 
	\begin{align}\label{weakformF2} \notag
		\skp{\us(t)}{\phi}_{\lpr{2}{\cO}} 
		=
		\skp{\uen}{\phi}_{\lpr{2}{\cO}}
		&+
		\int_{0}^{t \wedge T_{\sigma}} \skp{\fast}{\phi}_{H^{1}(\cO)} \ds
		\\ 
		&+
		\sum_{ k \in \Z} \int_{0}^{t\wedge T_{\sigma}}  \skp{\ssl}{\phi}_{\lpr{2}{\cO}} \dbl
	\end{align}
	for all $\phi \in \Hsobper{1}{\cO}$. Both formulations \eqref{weakformF1} and \eqref{weakformF2} hold for all $t \in [0,T]$. 
	
	Let us verify \eqref{Krylov:StochInt} for the stochastic integral in \eqref{weakformF2}.
	With \eqref{eq:ONB-1}-\eqref{eq:ONB-3} and the assumptions on the data \eqref{Lambda-Wachstumsbedingung} we get
	\begin{align*}
		\sum_{k \in \Z} \int_{\cO} \lal^{2} \abs{(\ue \gl)_{x}}^{2} \dx
		=  \frac{\lan^{2}}{L} \int_{\cO} (\uex)^{2} \dx 
		+
		\sum_{k = 1}^{\infty} \lal^{2}\frac{2}{L} \int_{\cO} (\uex)^{2} \dx
		+
		\sum_{k = 1}^{\infty}   \lal^{2} k^{2} \frac{8 \pi^{2}}{L^{3}} \int_{\cO} (\ue)^{2} \dx \, ,
	\end{align*}
	which implies \eqref{Krylov:StochInt} due to \eqref{est:EnEntropFischeGruen}.
	The proof for $(\ssl)_{x}$ uses similar arguments and will be omitted.
	
	For convenience, we will state the operators we work with once more.
	To guarantee well-posedness on $H$ and continuity of their Fr\'echet derivatives on $H \times H$, we use a cut-off function $\eta \in C^{\infty}(\R)$ 
	such that for an appropriate $\delta > 0$ 
	\begin{align}\label{def:CutoffFunc}
		\eta(x) =
		\begin{cases}
			\abs{x} &\qquad \text{for} \quad \abs{x} \ge \bar{c}_{\e}\sigma^{\frac{2}{p-2}}
			\\
			\in \R^{+}  &\qquad \text{for} \quad \abs{x} \in (\bar{c}_{\e}\sigma^{\frac{2}{p-2}} -\delta, \bar{c}_{\e}\sigma^{\frac{2}{p-2}})
			\\
			\bar{c}_{\e}\sigma^{\frac{2}{p-2}} - \delta &\qquad \text{for} \quad \abs{x} \le \bar{c}_{\e}\sigma^{\frac{2}{p-2}} - \delta
		\end{cases}
	\end{align}
	and 
	\begin{align}\label{est:CutoffFuncDer}
		\abs{\eta^{(s)}(x)} \le C(s) \delta^{-s} \quad, \, s\in (1,2) \, .
	\end{align}
	 Note that \eqref{est:ue-positive} implies $\eta(\ue) =\ue$, $\eta'(\ue) =1$, and $\eta''(\ue) = 0$ on $[0,T_{\s}]$. We define 
	\begin{align}\label{def:Energy-Operator}
		&E_{1}: u \mapsto \frac{1}{2}\int_{\cO} u^{2} \dx \, ,
		\\ \label{def:Potential-Operator}
		&E_{2}: u \mapsto \e \int_{\cO} \W (\eta(u)) \dx \, ,
		\\ \label{def:Entropy-Operator}
		&\Gmcas: u \mapsto \int_{\cO} G_{\alpha}(\eta(u)) \dx \, ,
	\end{align}
	where $G_{\alpha}$ is defined by 
	\begin{align}
		G_{\alpha}(u) = \frac{1}{\alpha(\alpha +1)} u^{\alpha+1} - \frac{1}{\alpha}u + \frac{1}{\alpha+1} 
		> 0 \, .
	\end{align} 
	For $E_{1}$, $E_{2}$, and $\Gmcas$, we compute the Fr\'echet derivatives 
	\begin{align} \label{FrechetDer-E1}
		DE_{1} (u)[v] 
		= 
		\int_{\cO} uv \dx
		\qquad  
		D^{2}E_{1} (u)[v,w] 
		= 
		\int_{\cO} vw \dx
	\end{align}
	\begin{align}\label{FrechetDer-E2} 
		DE_{2} (u)[v] 
		= 
		\e \int_{\cO} \W'(\eta(u)) \eta'(u) v \dx
	\end{align}
	\begin{align}\label{FrechetDer-E22} 
		D^{2}E_{2} (u)[v,w] 
		= 
		\e\int_{\cO} \W''(\eta(u)) (\eta'(u))^{2} vw + \W'(\eta(u)) \eta''(u) vw \dx
	\end{align}
	\begin{align}\label{FrechetDer-G}
		D\Gmcas (u)[v] 
		= 
		\int_{\cO} \Gas'(\eta(u))\eta'(u)v \dx
		= 
		\int_{\cO}
		\left(
		\frac{1}{\alpha}\eta(u)^{\alpha} - \frac{1}{\alpha}
		\right)
		\eta'(u)  v 
		\dx
	\end{align}
	\begin{align}\label{FrechetDer-G2}
		D^{2}\Gmcas (u)[v,w] 
		&= 
		\int_{\cO} \Gas''(\eta(u)) (\eta'(u))^{2} vw 
		+ 
		\Gas'(\eta(u)) \eta''(u) vw \dx \, .
	\end{align}
	Due to the cutoff $\eta$, the assumptions i) to iv) in Theorem \ref{theo:Ito/KrlyovFormel} are readily checked for the space $H$ and its dense subsets $V_{1}$ (in the case of $E_{1}$) and $V_{2}$ (in the case of $E_{2}$ and $\Gmcas$, respectively).
	Hence, we may choose the operators $E_{1}, E_{2}$ and $\Gmcas$ to use It\^o's formula w.r.t. to the weak formulation \eqref{weakformF1} in the first case and w.r.t. \eqref{weakformF2} in the other two cases. We end up with the following equations which hold for $t \in [0,\Tmax]$.
	\begin{align}\label{eq:ItoE1a}\notag
		\frac{1}{2}\int_{\cO} (\us(t))^{2}_ {x} \dx
		&=
		\frac{1}{2}\int_{\cO} (\uen)^{2}_{x} \dx
		+
		\sum_{ k \in \Z} \int_{0}^{t\wedge T_{\sigma}} \int_{\cO} (\us)_{x}   \lal (\us\gl)_{xx} 
		\dx \dbl 
		\\ \notag
		&+
		\int_{0}^{t\wedge T_{\sigma}} \int_{\cO} -(-(\us)^{2} \psx - C_{Strat} (\us)_{x}) (\us)_{xxx} \dx \ds
		\\ 
		&+
		\frac{1}{2} \int_{0}^{t\wedge T_{\sigma}}  \sum_{ k \in \Z} \lal^{2}\int_{\cO}(\us \gl)^{2}_{xx} \dx \ds \, 	
	\end{align}  
	and
	\begin{align}\label{eq:ItoE2a}\notag
		\e \int_{\cO} \W(\us(t)) \dx  
		&=
		\e \int_{\cO} \W(\eta(\uen))) \dx
		+
		\sum_{ k \in \Z} \int_{0}^{t\wedge T_{\sigma}} \int_{\cO} \e \W'(\us)  \lal (\us \gl)_{x} \dx \dbl
		\\ \notag
		&+
		\int_{0}^{t\wedge T_{\sigma}} \int_{\cO}
		(-(\us)^{2} \psx - C_{Strat} (\us)_{x}) (\e \W'(\us))_{x}
		\dx \ds
		\\
		&+
		\frac{1}{2} \int_{0}^{t\wedge T_{\sigma}}  \sum_{ k \in \Z} \lal^{2}\int_{\cO} \e \W''(\us)(\us \gl)^{2}_{x} \dx \ds \, .	
	\end{align} 	
	The entropy \eqref{def:Entropy-Operator} applied to \eqref{weakformF2} then gives
	\begin{align}\notag
		\label{eq:ItoGa1}
		\Gmcas (\us(t)) 
		=
		\Gmcas (\uen) 
		+
		&\sum_{ k \in \Z} \int_{0}^{t\wedge T_{\sigma}}\int_{\cO} 
		\left(
		\frac{1}{\alpha} (\us)^{\alpha} - \frac{1}{\alpha}
		\right)
		\lal(\us \gl)_{x} 
		\dx \dbl
		\\ \notag
		+
		&\int_{0}^{t\wedge T_{\sigma}} \int_{\cO} 
		(-(\us)^{2} \psx - C_{Strat} (\us)_{x})
		\left(\frac{1}{\alpha} \us^{\alpha} - \frac{1}{\alpha}
		\right)_{x}
		\dx \ds
		\\ 
		+
		&\frac{1}{2} \int_{0}^{t\wedge T_{\sigma}} \sum_{ k \in \Z} \lal^{2}\int_{\cO}\us^{\alpha-1}(\us \gl)^{2}_{x} \dx \ds \,.
	\end{align}
	
	Thus, combining \eqref{eq:ItoE1a}, \eqref{eq:ItoE2a}, and \eqref{eq:ItoGa1} and adding all terms with a good sign to the left hand side, we obtain equation \eqref{eq:CombinedEnEntr}.
	
	\section{Predictability}\label{predictability}
	In this section we prove predictability of the processes $\xiast$ and $\fast$, cf. \eqref{RieszProjH2} and \eqref{RieszProjH1}, as well as of $\ssl$ and $(\ssl)_{x}$.
	Let us first show that with $u$ also $\ux$ is predictable, where we skipped the index $\e$. By means of the Skorokhod-embedding theorem we will further assume that we are operating on a probability space $([0,1],\B([0,1]),\lambda) = (\Omega,\Sigma, \mathbb{P})$, where $\lambda$ denotes Lebesgue measure and $\B([0,1])$ the Borel $\sigma$-algebra on $[0,1]$. We write  $\mathcal{P}(\Omega \times [0,T])$ for the predictable $\sigma$-algebra on $\Omega \times [0,T]$.
	By a convolution argument and an application of Theorem 15.1 in \cite{Kechris1995} and Sobolev embedding, we obtain  
	\begin{lem1}\label{predictability ux}
		Let $u \in \lpr{1}{\Omega; \lpr{1}{[0,T];\lpr{s_{1}}{\cO}}}$ be predictable as a mapping \newline $u \colon \Omega \times [0,T] \rightarrow \lpr{s_{1}}{\cO}$ and let $\ux \in \lpr{1}{\Omega; \lpr{1}{[0,T]; \lpr{s_{2}}{\cO}}}$. Then, $\ux \colon \Omega \times [0,T] \rightarrow \lpr{s_{2}}{\cO}$ is predictable. 
	\end{lem1}
	Combining \eqref{est:EnEntropFischeGruen} and \eqref{lem:UniformEstimate-u^2ph2} with Lemma \ref{predictability ux}, we get
	\begin{kor1}
		$\ssl$ and $(\ssl)_{x}$ are predictable as mappings into $\lpr{2}{\cO}$.
	\end{kor1}
	To establish the predictability of $p_{x}$, we first propose the following integrability result. 
	\begin{lem1}\label{higher integrability of px}
		We have that $p_{x}$ is an element of the space $\lpr{\eta}{\Omega \times [0,T]; \lpr{\eta}{\cO}}$ for $\eta < \frac{6p}{2+3p}$.
	\end{lem1} 
	It is based on standard estimates combined with Lemma \ref{lem:Gagliardo-Nirenberg-est} and Corollary \ref{kor:higher_reg_p}.
	\begin{lem1}\label{lem:Gagliardo-Nirenberg-est}
		Let $q \in [2,6)$ and \newline $v \in \lpr{\frac{2(q+2)}{6-q}}{\Omega;\lpr{\infty}{[0,T];\lpr{2}{\cO}}} \cap \lpr{2}{\Omega;\lpr{2}{[0,T];\Hsob{1}{\cO}}}$. Then we have
		\begin{align}\notag
			\Norm{v}_{\lpr{q}{\Omega;\lpr{q}{[0,T];\lpr{q}{\cO}}}} 
			&\le
			C T^{\frac{6-q}{4q}}
			\left(
			\Norm{v_{x}}^{\frac{q-2}{2q}}_{\lpr{2}{\Omega;\lpr{2}{[0,T];\lpr{2}{\cO}}}} \cdot \Norm{v}^{\frac{q+2}{2q}}_{\lpr{\frac{2(q+2)}{6-q}}{\Omega;\lpr{\infty}{[0,T];\lpr{2}{\cO} } } } 
			\right) 
			\\
			&\quad + C\Norm{v}^{2}_{\lpr{2}{\Omega;\lpr{2}{[0,T];\lpr{2}{\cO}}}} \,.
		\end{align}
	\end{lem1}
	\begin{kor1}\label{kor:higher_reg_p}
		For $\eta < 3p$ we have
		\begin{align}\label{kor:higher_reg_p_h1}
			\int_{\Omega} \int_{0}^{T}\int_{\cO} u^{-\eta} \dx \ds \dif \mathbb{P}  \le C
		\end{align} 
		and in particular 
		\begin{align}\label{kor:higher_reg_p_h2}
			\int_{\Omega} \int_{0}^{T}\int_{\cO} u^{-2p-2} \dx \ds \dif \mathbb{P}  \le C \, .
		\end{align} 
	\end{kor1}
	\proof
	By the energy-estimate in Theorem \ref{theo:existence-stfe-stratonovich}, we can apply Lemma \ref{lem:Gagliardo-Nirenberg-est} to $u^{-\frac{p}{2}}$ and infer that for arbitrary $q \in (0,6)$
	\begin{align}\label{kor:hihgher_reg_p_h3}
		\notag
		&\quad \int_{\Omega} \int_{0}^{T}\int_{\cO} u^{-\frac{p}{2}\cdot q} \dx \dt \dif \mathbb{P}
		\\ \notag
		&\le
		C T^{\frac{6-q}{4}}
		\left(
		\Norm{(u^{-\frac{p}{2}})_{x}}^{\frac{q-2}{2}}_{\lpr{2}{\Omega;\lpr{2}{[0,T];\lpr{2}{\cO}}} }
		\cdot
		\Norm{u^{-\frac{p}{2}}}^{\frac{q+2}{2}}_{\lpr{\frac{2(q+2)}{6-q}}{\Omega;\lpr{\infty}{[0,T];\lpr{2}{\cO} } } }
		\right)
		\\
		&\quad + C\Norm{u^{-\frac{p}{2}}}^{2q}_{\lpr{2}{\Omega;\lpr{2}{[0,T];\lpr{2}{\cO}}}} \,.
	\end{align}
	The second factor in the first term and the second term are controlled by the energy estimate \eqref{est:EnEntropFischeGruen} in Theorem \ref{theo:existence-stfe-stratonovich}. Moreover, we have
	\begin{align}
		\Norm{u^{-\frac{p}{2}}}_{\lpr{2}{\Omega;\lpr{2}{[0,T];\Hsob{1}{\cO}}}} 
		\le C +
		C\int_{\Omega} \int_{0}^{T}\int_{\cO} u^{-p-2}\ux^{2} \dx \ds \dif \mathbb{P} 
		\le C \, ,
	\end{align}  
	where the last estimate follows with Lemma \ref{lem:UniformEstimate-u^2p}.
	Choosing $q = \frac{4p+4}{p}$ in \eqref{kor:hihgher_reg_p_h3}, we get \eqref{kor:higher_reg_p_h2}.
	\proofend
	With the results established above, we can prove
	\begin{lem1}\label{lem:p_x_predictable}
		For $\eta < \frac{6p}{2+3p}$, $p_{x}$ is predictable as a mapping from $\Omega \times [0,T] $ to $\lpr{\eta}{\cO}$.
	\end{lem1} 
	\proof
	From \eqref{lem:UniformEstimate-u^2ph2} in Lemma \ref{lem:UniformEstimate-u^2p} we have $\uxx \in \lpr{2}{\Omega;\lpr{2}{[0,T];\lpr{2}{\cO}}}$. 
	Since by \eqref{kor:higher_reg_p_h1} also $u^{-p} \in \lpr{2}{\Omega; \lpr{2}{[0,T]; \lpr{2}{\cO}}}$ holds, we see that $p: \Omega \times [0,T] \rightarrow \lpr{2}{\cO}$ is 
	%(up to equivalence classes) 
	well defined. Using Lemma \ref{predictability ux} twice, the predictability of $\uxx$ follows. Moreover, $u >0$ implies that $u^{-p}$ is predictable, as it is the composition of a predictable and a continuous mapping. Thus, also $p$ is predictable. Now we may infer as in Lemma \ref{predictability ux} that $p_{x}$ is predictable as a mapping into $\lpr{\eta}{\cO}$, $\eta < \frac{6p}{2+ 3p}$, where we used the $\lpr{\eta}{\Omega \times [0,T]; \lpr{\eta}{\cO}}$-integrability of $p_{x}$ established in Lemma \ref{higher integrability of px}.
	\proofend
	\begin{bem1}\label{rem:predict}
		The standard product of functions in $\Hsob{1}{\cO}$ and $\lpr{\eta}{\cO}$ conserves $\lpr{\eta}{\cO}$-integrability. 
		Moreover,  $\B(\Hsob{1}{\cO}) \subset \B(\lpr{\eta}{\cO})$. 
		Thus, the mapping $\tau:\Hsob{1}{\cO} \times \lpr{\eta}{\cO} \rightarrow \lpr{\eta}{\cO}$ with $\tau(v,w) := v\cdot w$ is well defined and continuous. 
		Using similar arguments as in the finite dimensional case, cf. for example \cite{Klenke2014eng} Chapter 1.4, one finds the mapping $h: \Omega \times [0,T] \rightarrow \Hsob{1}{\cO} \times \lpr{\eta}{\cO}$, $(\omega,t) \mapsto (v(\omega,t),w(\omega,t))$, to be predictable, if $v$ and $w$ are. From Lemma \ref{predictability ux} and the continuity of $u$, we infer that $u$ is  predictable as a mapping from $\Omega \times [0,T]$ into $\Hsob{1}{\cO}$. Using Lemma \ref{lem:p_x_predictable} as well, the predictability of $u p_{x} = \tau \circ h(u,p_{x})$ follows, as it is the composition of a predictable and a continuous mapping.
		Iterating this argument, we obtain the predictability of $u^{2}p_{x}$ as a mapping from $\Omega \times [0,T]$ into $\lpr{\eta}{\cO}$.
	\end{bem1}
%	Finally, we show that 
	By a convolution argument combined with Remark \ref{rem:predict}, we obtain the following result. 
	\begin{lem1}
		$up_{x}$ and $u^{2}p_{x}$ are predictable as mappings from $\Omega \times [0,T]$ into $\lpr{2}{\cO}$. 
	\end{lem1}

	Now, again by continuity the predictablility of $\xiast$ and $\fast$ follows.
	\begin{lem1}
		The mappings $\xiast$ and $\fast$ are predictable as mappings from $\Omega \times [0,T]$ to $\Hsobper{2}{\cO}$ and $\Hsobper{1}{\cO}$, respectively.
	\end{lem1}
	\proof
	Let us set $X := \lpr{2}{\Omega;\lpr{2}{[0,T]; \lpr{2}{\cO}}}$ and $Y := \lpr{2}{\Omega; \lpr{2}{[0,T];(\Hsobper{2}{\cO})'}}$.
	We first consider the mapping $\eta:X \rightarrow Y$, with $\eta(u) = \zeta$, where 
	\begin{align}
		\zeta(v) := \Ew{\int_{0}^{T} \int_{\cO} u_{xx}  v \dx \dt }
	\end{align}
	for $v \in Z := \lpr{2}{\Omega; \lpr{2}{[0,T];\Hsobper{2}{\cO}}}$. Due to the periodic boundary conditions, by integration by parts one finds that $\zeta$ and consequently also $\eta$ is well defined. We now show that $\eta$ is also continuous. Let $u,w$ be elements in $X$ such that for $\delta > 0$, $\Norm{u-w}_{X} < \delta$ holds. Then, by means of Hölder's inequality and  integration by parts, we may estimate
	\begin{align}\notag
		&\underset{\Norm{v}_{Z }=1}{\sup} 
		\left \vert
		\Ew{\int_{0}^{T} \int_{\cO}  (u-w)_{xx} v \dx \dt }
		\right \vert
		\\ 
		&\le
		\Ew{\int_{0}^{T} \int_{\cO}  (u-w)^2 \dx \dt }^{\frac{1}{2}}
		\underset{\Norm{v}_{Z} =1 }{\sup}
		\Ew{\int_{0}^{T} \int_{\cO} (v_{xx})^{2} \dx \dt }^{\frac{1}{2}}
		\le C \delta \, ,
	\end{align} 
	which implies that $\eta$ is continuous. Since $\xiast$ is obtained from $\xi$ via the isometric isomorphism $J:Z \rightarrow Y$ from Riesz' representation theorem, cf. \eqref{RieszProjH2}, $\xiast$ is predictable as a mapping from $\Omega \times [0,T]$ to $\Hsobper{2}{\cO}$, as it is a continuous transformation of the predictable mapping $-u^{2} p_{x} - \Cstr \ux$:
	\begin{align}
		\xiast = (J^{-1}\circ \eta)(-u^{2} p_{x} - \Cstr \ux) \, .
	\end{align}
	For the predictability of $\fast$ one can argue in the same manner, using the corresponding spaces w.r.t. to $\Hsobper{1}{\cO}$.  
	\proofend

\noindent
{\bf Acknowledgment.} L.K. has been supported by the
Research-Training-Group 2339 “Interfaces, Complex Structures, and Singular Limits” of German
Research Foundation (DFG). G.G. gratefully acknowledges in addition the support of DFG through the project ``Free boundary
propagation and noise: analysis and numerics of stochastic degenerate parabolic equations’.

\end{document}